\documentclass[11pt]{report} 
\usepackage{amsmath,
            amsfonts,
            amssymb,
            amstext,
            latexsym,
            times,      
            verbatim,   
            }

\def \dfll {\leaders \hbox to 1em {\hss.\hss}\hfill}
\newtheorem{theorem}{Theorem}
\newtheorem{propos}[theorem]{Proposition}
\newtheorem{cor}[theorem]{Corollary}

\newtheorem{definition}{Definition}
\newtheorem{qes}{Question}
\newtheorem{lemma}{Lemma}
\def\blm{\begin{lemma}}
\def\elm{\end{lemma}}
\def\bdf{\begin{defi}}
\def\edf{\end{defi}}
\def\btm{\begin{theorem}}
\def\etm{\end{theorem}}
\def\bpp{\begin{propos}}
\def\epp{\end{propos}}
\def\bQ {\begin{qes}}
\def\eQ {\end{qes}}
\def\btm{\begin{theorem}}
\def\etm{\end{theorem}}
\def\ben{\begin{enumerate}}
\def\een{\end{enumerate}}

\def \Mu {{\sf U}}

\def\m {{\sf m}}

\newcommand{\figref}[1]{(\ref{fig-#1})}

\def\Cob {\mbox{\boldmath${\cal C}ob$\ \unboldmath}}
\def\Alg {\mbox{\boldmath${\cal A}lg$\ \unboldmath}}
\def\Tgl {\mbox{\boldmath${\cal T}\!gl$\ \unboldmath}}

\def\ep{\ \hfill{\rule {2.5mm}{2.5mm}}\smallskip}

\newcommand{\lbl}[1]{\label{#1}}

\newcommand{\picbox}[1]{\begin{equation}\label{fig-#1}
\epsfbox{#1.eps}\end{equation}}

\def \isto {\widetilde{\longrightarrow}}

\newcommand{\TO}[2]{\stackrel {\mbox{#1}}{\hbox to #2pt{\rightarrowfill}}}

\def\thrafill{$\mathsurround=0pt \mathord- \mkern-6mu 
\cleaders\hbox{$\mkern-2mu
\mathord- \mkern-2mu$}\hfill \mkern-6mu\mathord\twoheadrightarrow$}
\newcommand {\onto} [1]{\hbox to #1pt{\thrafill}}

\usepackage{amsfonts,
            amssymb}

\input{epsf.sty}

\newcommand{\head}[1]{\medskip \begin{center}{\large \sc #1}\end{center}}                     

\newcommand{\Z}{{\mathbb Z}}

\newcommand{\R}{{\mathbb R}}

\newcommand{\N}{{\mathbb N}}

\newcommand{\posit}[2]
{\raise -1.4ex\hbox{${\textstyle #1}\atop {\stackrel{\uparrow}{#2}}$}}

\def\Ii {\mbox{\raise .4 ex\hbox{$\int$}$\!\! I$}}

\begin{document} 

\vspace*{1cm}

\begin{center}

\section*{Towards an Algebraic Characterization of 3-dimensional Cobordisms}

\bigskip

\medskip

{\large Thomas Kerler}\\
\medskip

April 2001
 \vspace*{1.5cm}

(To appear in {\em Contemp. Math.})

\bigskip

\end{center}

{\small \noindent{\bf Abstract :}
The goal of this paper is to find a close to isomorphic presentation of 
3-manifolds in terms of Hopf algebraic expressions. To this end 
we define and compare three different braided tensor categories that arise
naturally in the study of Hopf algebras and 3-dimensional topology. 
The first is the category \Cob of connected surfaces with 
one boundary component and 3-dimensional relative cobordisms, the second
is a category \Tgl of tangles with relations, and the third is 
a natural algebraic category \Alg freely generated by a Hopf algebra 
object. From previous work we know that {\Tgl} and 
\Cob are equivalent. We use this fact and the idea of Heegaard splittings 
to construct a surjective 
functor  from \Alg  onto  {\Cob}. We also find a 
map that associates to the generators of the mapping class group in 
\Cob preimages in \Alg. The single block 
relations in the mapping class group are verified for these expressions. 
We propose to find a version of  \Alg with possibly additional
relations to obatin isomorphic algebraic presentations of the mapping class groups and
eventually of \Cob. 
\footnote{Lecture presented at the 958th Meeting of the AMS, 
San Francisco State University,  Session on Diagrammatic Morphisms in
Algebra, Category Theory, and Topology.

2000 Mathematics Subject Classification: 
Primary 57R56; Secondary 14D20, 16W30, 17B37, 18D35, 57M27.}

}
\medskip

\begin{center}
{\large \sc Contents}\
\smallskip

\parbox[t]{11cm}{
1. Introduction\dotfill\pageref{S1}

2. The category \Cob\dotfill\pageref{S2}

3. Presentation of \Cob by \Tgl\dotfill\pageref{S3}

4. The category \Alg\dotfill\pageref{S4}

5. A functor from \Alg to \Cob\dotfill\pageref{S5}

6. From generators of \Cob to generators of \Alg\dotfill\pageref{S6}

}
\end{center}

\head{1. Introduction}\lbl{S1}

  For some time is has been a puzzling question whether the 
appearance of Hopf algebras in the world of quantum invariants of
3-manifolds as in \cite{ResTur91}
is a lucky coincidence that makes computations work
or if these structures arise in more fundamental ways out of 
3-dimensional topology. 

 It was soon understood that such algebraic structures 
are in fact inherent in 
the category of cobordisms  \Cob between connected surfaces with
one boundary component. Specifically, the torus with one hole as an 
object in \Cob was discovered to admit the
structure of a braided Hopf algebra in the sense of \cite{Maj93}. 

  The three-dimensional pictures of the cobordisms representing the 
products and coproducts have been found by Crane and Yetter and proven to
satisfy the relevant axioms, see \cite{Yet97} and \cite{CY99}. 
The same picture emerged independently in investigations by the author
via the route of tangle presentations,
see \cite{Kerun}. An interpretation of braided Hopf algebra structures in 
linear abelian braided tensor categories was found by Lyubashenko in 
\cite{Lub95}. This combined with the equivalence of  cobordisms categories 
with certain tangle categories given in \cite{Ker99} leads to the three 
dimensional interpretation. The cobordisms are easily worked out to be
the same as the ones in  \cite{Yet97} and \cite{CY99}, see \cite{Kerun}.
The tangle pictures are also used as examples for integrals in braided
categories in \cite{BKLT00}. 

The purpose of this note is to prove a theorem that was already stated 
in \cite{Ker96}, which not only asserts the existence of cobordisms representing
structure morphisms such as products and coproducts of a Hopf algebra but also that
these generate the entire cobordism category. In more formal terms we will define
an algebraic category \Alg via generators and relations representing the
axioms for a braided Hopf algebra in a braided tensor category. The existence of
the special cobordisms is implied by the existence of a functor and their generating
property by  surjectivity of this functor. 

\btm\label{thm-main} Let \Alg be the braided tensor category freely generated by a 
Hopf algebra object as defined in Section~4, and \Cob the cobordism category defined
in Section~2.
There exists a surjective functor of braided tensor categories 
$$
{\mathfrak G}\;:\;\;\Alg\;\longrightarrow\;\Cob\;.
$$ 
\etm

The proof of this theorem as we present it here uses and demonstrates several
techniques in graphical and diagrammatic categorical calculations. Particularly,
we give a description of \Alg as a category of directed trivalent graphs with 
crossings and special types of endpoints modulo relations. Also, in Section~3
we recall how the cobordism category \Cob can be presented in terms of a
category of tangles modulo relations \Tgl. The constructions of functors and
assignments in Section~5 and  Section~6 are done entirely in these graphical
languages. 

In  Sections~6.2 and 6.3 we also discuss the problem of modifying \Alg so
that ${\mathfrak G}$ becomes an isomorphism of categories. This means we would
have to find an assignment of generators of \Cob to generators of \Alg such that
the relations in \Cob are also respected in \Alg. For the genus one relations this
is in fact true but for higher genera it is likely that we will have to impose more
relations on \Alg. 
We summarize next the observations we will to make on
this question.

\btm\label{thm-more} 
\ben
\item Let ${\Gamma}_{1,1}^*$ be the central extension of the mapping class group
of the torus with one hole.
Let $A\in \Alg$ be the generating Hopf algebra object. 
 Then there exists a homomorphism 
${\mathfrak W}^{[1]}:\widehat{\Gamma}_{1,1}\to Aut(A)$, such
that the composite
$$
{\Gamma}_{1,1}^*
\;\stackrel{{\mathfrak W}^{[1]}}{-\!\!\!-\!\!\!-\!\!\!-\!\!\!\longrightarrow}\;
Aut(A)
\; \stackrel{{\mathfrak G}}{-\!\!\!-\!\!\!-\!\!\!-\!\!\!\longrightarrow}\;
\widehat{\Gamma}_{1,1}\subset \Cob
$$
is the identity. Here $Aut(A)$ denotes the invertible morphisms $A\to A$ in 
\Alg, and we use  that ${\Gamma}_{1,1}^*$ is identical with
the group of invertible cobordisms on the surface of genus one. 
\item There is a natural set $Gen[\Cob]$ of generators of \Cob and an assignment 
${\cal W}:Gen[\Cob]\to \Alg$ such that ${\mathfrak G}\circ {\cal W}$ is the 
identity on $Gen[\Cob]\subset\Cob$. 
\item Suppose it is possible to find additional relations on \Alg, yielding a 
subquotient category $\overline{\Alg}$ such that ${\mathfrak G}$ factors into
$\overline{\Alg}$ and ${\cal W}$ extends to a functor ${\mathfrak W}$
on $\overline{\Alg}$. 

Then ${\mathfrak W}$ is a two-sided inverse of ${\mathfrak G}$, and hence 
$\overline{\Alg}\cong \Cob$. 
\een
\etm

The challenging question now is to find the additional relation that
define $\overline{\Alg}$. This would imply a characterization of three
dimensional topology in purely algebraic terms!

\paragraph{Acknowledgements:} I would like to thank the organizers
of the AMS session,  David Radford,  Fernando Souza, and
David Yetter, for the opportunity to speak about this paper. I also
would like to thank Volodya Lyubashenko for discussion and Jim Stasheff
for encouraging me to write up these results following \cite{Kerun}.

\head{2. The category \Cob :}\lbl{S2}

\paragraph{2.1 Category of 2-framed, relative cobordisms} 
For every integer $g\geq 0$ construct a model surface $\Sigma_g$ of genus $g$
with  $\partial \Sigma_g=S^1$. This can be done as follows. Pick a 
disc $D^2$ and cut out $2g$
small discs along a diameter of $D^2$. Along the new boundary components glue 
in $g$ cylinders such that the ends of
 each cylinder are glued to two consecutive holes. The result
is depicted below. 
\picbox{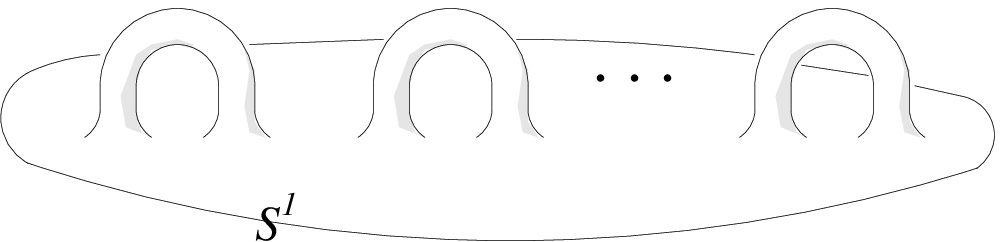}

For two surfaces $\Sigma_h$ and $\Sigma_g$ let   
$\Sigma_{[h,g]}=-\Sigma_h\cup_{S^1\times 0}S^1
\times[0,1]\cup_{S^1\times 1}\Sigma_g$
be the closed oriented surfaces of genus $g+h$ obtained by gluing the original  
surfaces together along their boundaries with a cylinder inserted. A 
{\em relative cobordism} is a compact oriented manifold $M$, together with a
homeomorphism $\psi:\partial M\isto \Sigma_{[h,g]}$.

We consider two pairs of data $(M,\psi)$ and $(M',\psi')$ as equivalent if there
is a homeomorphism $h:M\isto M'$ such that $\psi'\circ h=\psi$ on $\partial M$.
We write the equivalence class $[M,\psi]$, with minor abuse of notation, in the
morphism form:
$$
M\,:\;\Sigma_h\,\longrightarrow\,\Sigma_g
$$
It is not hard to see (e.g., \cite{KerLub00}) that $[M,\psi]$ does not 
change under isotopies of $\psi$.

For two cobordisms $M:\Sigma_g\to\Sigma_h$ and $N:\Sigma_h\to\Sigma_k$ we 
define a composite $N\circ M:\Sigma_g\to\Sigma_k$ by gluing two cobordisms
together along the common boundary piece $\Sigma_h$ using the coordinate
maps on $\partial M\isto \Sigma_{[h,g]}$ and $\partial N\isto \Sigma_{[g,k]}$.
The cylindrical part of length  2 is then monotonously shrunk to a cylinder of
height 1. It is easy to see that the resulting homeomorphism class does not
depend on the choices of $M$ and $N$ in their classes.  

Thus we obtain a category, denoted $\Cob_0$, 
which has the surfaces $\Sigma_g$ as objects
and the classes $[M,\psi]$ as morphisms.  

For a relative cobordism $M:\Sigma_h\to\Sigma_g$ denote by $M_0$ the 
cobordism obtained by gluing in a full cylinder $D^2\times [0,1]$ along
the $S^1\times [0,1]$ part of the boundary so that $M_0$ is a cobordism 
between closed surfaces $\overline{\Sigma_g}=\Sigma_g\cup D^2$. 
Consider an unknotted embedding of $\overline{\Sigma_g}\times [0,\epsilon]$
into $S^3$. 
A framing on $S^3$ hence induces a standard framing on the collar
$\overline{\Sigma_g}\times [0,\epsilon]$. 

In addition to the topological structure from above we 
 consider now also manifolds with 2-framings, i.e., isotopy classes
of trivializations $TM\oplus TM \isto \R^6\otimes M$, and assume standard
trivializations on standard handle bodies bounding the surfaces. We restrict
to those that are compatible with the standard 2-framings on the collars
 $\overline{\Sigma_g}\times [0,\epsilon]$. As a result the gluing 
composition operation extends to the 2-framed cobordisms. 

It is standard knowledge, see for example \cite{Ker99},
 that this information is equivalent to
the signature of a bounding 4-manifold. We thus obtain an exact sequence 
in the sense of group theory: 
$$
1\,\longrightarrow\,\Z\,\longrightarrow\,\Cob\, \longrightarrow\,\Cob_0\,
\longrightarrow\,1\qquad . 
$$

\paragraph{2.2 Braided tensor category:}

The category $\Cob_0$ has a natural tensor product. It is given on the 
objects by $\Sigma_g\otimes\Sigma_h=\Sigma_{g+h}$. In order to describe $\otimes$
on the morphisms we choose a
disc $P=D^2-(D^2\sqcup S^2)$ 
with two holes. The   two surfaces $\Sigma_g$ and $\Sigma_h$ 
are sewn them into the two holes of $P$ 
such that their handles are aligned as depicted below. Upto isotopy there is
then a unique homeomorphism $\lambda_{g,h}:
\Sigma_g\cup_{S^1}P \cup_{S^1}\Sigma_h\,\isto\,\Sigma_{g+h}$ which maps 
the corresponding  handles in order onto each other.  

\picbox{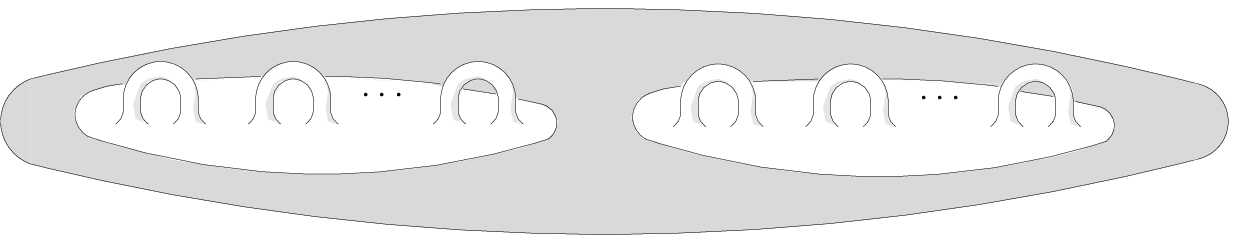}

For two cobordisms $M:\Sigma_h\to\Sigma_g$ and 
$L:\Sigma_p\to\Sigma_q$ the tensor product is obtained by gluing the 
cobordisms into $P\times [0,1]$ and using the $\lambda_{g,h}$ to adjust the
boundary identifications so that $M\otimes N:\Sigma_{h+p}\to\Sigma_{g+q}$. 

Since $\lambda_{a,b+c}(1\otimes \lambda_{b,c})$ is isotopic to 
 $\lambda_{a+b,c}( \lambda_{a,b}\otimes 1)$ this product is strictly 
associative. In \cite{KerLub00} we describe a procedure by which this 
tensor product lifts to the 2-framing structure so that we still have 
$(M\otimes L)\circ (N\otimes K)=(M\circ N)\otimes (L\circ K)$ for 
compatible cobordisms. 

 Finally, we obtain a family of isomorphisms 
$c_{g,h}:\Sigma_g\otimes\Sigma_h\isto\Sigma_h\otimes\Sigma_g$ from the cylinder
$P\otimes [0,1]$ by twisting the two ends by $\pi$ relative to each other
so that opposite holes in $P$ are connected to each other by boundary 
cylinders. The morphism $c_{g,h}$ is then obtained by gluing the cylinders
$\Sigma_g\times [0,1]$ and $\Sigma_h\times [0,1]$ into this twisted
version of $P$. See \cite{KerLub00} again for more details. 
By construction $c_{g,h}$ is a natural isomorphism from $\otimes T $
to $\otimes $ with all properties of a braiding.
Here $T:\Cob\times\Cob\to\Cob\times\Cob$ is the transposition
and $\otimes:\Cob\times\Cob\to\Cob$ as defined. The remarks so far 
are summarized as follows:
\btm
\Cob is a braided tensor category with objects $\Sigma_g=\Sigma_1^{\otimes g}$. 
\etm

\paragraph{2.3 Generators of \Cob :}

We will describe here the set of generators of \Cob that comes from 
Heegaard splittings of cobordisms. 

Let $H^+_0:\Sigma_0\to\Sigma_1$ be the cobordism, obtained by attaching 
a {\em full} handle $D^2\times [0,1]$ to a {\em thickening} 
$D^2\times [0,\epsilon]$  of the disc  from 2.1 along the holes in 
$D^2\times 0$. The boundary identification is such that the restriction of
the gluing construction to $S^1\times [0,1]$ is precisely the construction of
$\Sigma_1$ described in Section~2.1. 
Hence $H^+\cong D^2\times S^1$ is itself a full torus. 
We introduce the  cobordisms
\begin{equation}\label{eq-handlecob}
H^+_{g,k}:\,\Sigma_g\,\to\,\Sigma_{g+k}\;\qquad \mbox{with}\qquad
H^+_{g,k}\,= \,id_{\Sigma_g}\otimes \underbrace{H^+_1\otimes\ldots\otimes H^+_1}_k\;.
\end{equation}
They describe the {\em addition} 
of $k$ 1-handles to a surface of genus $g$.  
 
Conversely, we have a cobordisms $H^-_0:\Sigma_1\to\Sigma_0$ where we glue
a thickened disc to $\Sigma_1\times [0,1]$ along the longitude of 
$\Sigma_1\times 1$ so that $H^-_0\circ H^+_0=id_{\Sigma_0}\cong D^3$.  
 The cobordisms $H^-_{g,k}:\Sigma_{g+k}\to\Sigma_g$ are 
obtained analogously.

The second type of morphisms arise from the mapping class groups.
Consider the group $Homeo^+(\Sigma_g)$ 
of orientation preserving homeomorphims of $\Sigma_g$ to itself,
which leave the boundary pointwise fixed. The mapping class group of
$\Sigma_g$ is thus the group of path connected components, that is 
$\Gamma_{g,1}=\pi_0(Homeo^+(\Sigma_g))$. 

To an element $\psi\in Homeo^+(\Sigma_g)$ we assign a cobordism 
$I_{\psi}$ as follows. The representative cobordism is given by the cylinder
$\Sigma_g\times [0,1]$. The boundary identification with 
$\Sigma_{[g,g]}$ is the canonical map (identity) on $\Sigma_g\times 0$ 
and it is given by $\psi$ on $\Sigma_g\times 1$. The cobordism $I_{\psi}$
only depends on the isotopy class of $\psi$ in $\Gamma_{g,1}$, which
we abusively denote by the same letter. Now if $Aut(\Sigma)$ denotes the
group of invertible cobordisms on $\Sigma$ in \Cob we have the following
result:

\btm[\cite{KerLub00}]
The following map is an isomorphims of groups.
$$
\Gamma_{g,1}^*\,\longrightarrow\, Aut(\Sigma_g)\,:\;\;\psi\,\mapsto\, I_{\psi}
$$
\etm 
Here $\Gamma_{g,1}^*$ is the central extension of $\Gamma_{g,1}$, which 
carries the corresponding framing information.  

A Heegaard splitting of a cobordism is given now as follows:
\btm
Every cobordism $M:\Sigma_h\to\Sigma_g$ in \Cob is given as a composite
$$
M\;=\;H^-_{g,N-g}\circ I_{\psi}\circ H^+_{h,N-h}
$$
for some $N\geq max(g,h)$ and some $\psi\in\Gamma_{N,1}^*$. 
\etm

{\em Proof:} Consider the space of  differentiable functions $f:M\to [0,1]$ 
such that $f^{-1}(\Sigma_h)=\{0\}$, $f^{-1}(\Sigma_g)=\{1\}$, and on the 
cylindrical piece $\cong S^1\times [0,1]\subset\partial M$ $f$ coincides with
the canonical projection. Assume some metric on $M$.
By standard arguments from differential topology
we can assume that $f$ is a Morse function and has singularities only
of index 1 or 2. Further more we can assume that the critical values of index 2
all lie in $(\frac 12,1]$ and the critical values  of index 1 in 
$[0, \frac 12)$ and $\frac 12 $ is a regular value. 
With $\Sigma_{int}=f^{-1}(\frac 12)$  we have cobordisms 
$A=f^{-1}([0,\frac 12])$ and $B=f^{-1}([\frac 12,1])$ with $M=B\circ A$. 
The gradient flow of $f$ identifies $\Sigma_h$ with $2l$ discs removed
with a submanifold $\phi:\Sigma_h-2lD^2\hookrightarrow\Sigma_{int}$.
Here $l=N-h$ is the number of index 1 critical points and the 
locations of the discs are given by their critical manifolds. We can 
always from a map $\tau\in Homeo^+(\Sigma_h)$ isotopic to the identity,
which maps these discs into the standard position of the holes for
the next $l$ handle attachments for the construction of $\Sigma_{h+l}$.
There is up to isotopy a unique map that extends 
$\phi\circ\tau^{-1}:\Sigma_{g}-2lD^2\hookrightarrow \Sigma_{int}$ to a 
map $\hat\phi: \Sigma_{g+l}\isto \Sigma_{int}$ over the additional glued in 
cylinders. Using the gradient flow and its behavior around the singularities
the class of cobordism $A$ is given by  $I_{\hat\phi}\circ H^+_{h,l}$. 
By an analogous procedure we find $B=H^+_{g,k} \circ  I_{\hat\phi'}$,
which implies the claim if we set $\psi=\hat\phi'\circ \hat\phi$. 
The framing of $M$ is adjusted by choosing the appropriate extension class
in $\Gamma_{N,1}^*$. 

\ep

  The set of generators of \Cob can be broken down even further using 
the special generators of the mapping class groups. It is a well known
fact that $\Gamma_{g,1}$ is generated by a finite set of Dehn twists.
They are denoted by capital letters $A_j$, $B_j$ and $C_j$, for Dehn
twists along the curves depicted in \figref{mcg-gen} labeled by 
the corresponding lower case letters.      
\picbox{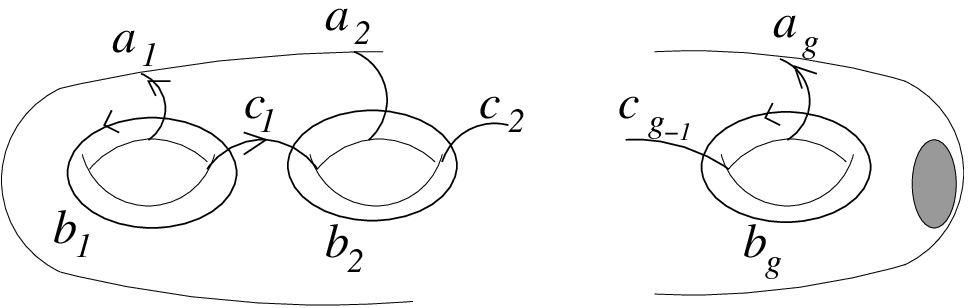}      

A slightly more convenient set of generators is given by the set 
$\{A_j,S_j, D_j\}$  where
\begin{equation}\label{eq-mcggen}
D_j\,=\,A^{-1}_jA^{-1}_{j+1}C_j\qquad\mbox{and}\qquad
S_j\,=\,A_jB_jA_j\;\qquad \;\; \mbox{for \ } j=1,\ldots,g. 
\end{equation}
In terms of cobordisms we can write 
$I_{A_j}=id_{\Sigma_{j-1}}\otimes I_{A_1} \otimes id_{\Sigma_{g-j}}\,$,
where $A_1\in \Gamma_{1,1}$. Similar formulae exists for $I_{S_j}$ 
and $I_{D_j}\,$. We imply here some specific representative in 
$\Gamma_{1,1}^*$. We also introduce the cobordism $Z:0\to 0$, which is topologically
the identity cylinder over $S^2$ but has framing changed by one. 
We find the following.
\begin{cor}\label{cor-gen}
As a tensor category \Cob is generated by the cobordisms 
$H^+_0$, $H^-_0$, $I_{A_1}^{\pm 1}$, $I_{D_1}^{\pm 1}$, $I_{S_1}^{\pm 1}$, and $Z^{\pm 1}$. 
\end{cor}

\head{3. Presentation of \Cob by \Tgl}\lbl{S3}

We recall a variant of the tangle presentation of \Cob given in 
\cite{Ker99}. 

\paragraph{3.1 Admissible  tangles and moves:}  
First we define the
category \Tgl. Its objects are non-negative integers. A morphism
$T:k\to l$ is obtained from a framed tangle  
in    $\R^2\times [0,1]$ with $2k$ end points 
$1^+, 1^-,\ldots, k^+, k^-$ at the top line $\R_x\times 1$ and $2l$ 
end points $1^+, 1^-,\ldots, l^+, l^-$ at the bottom line $\R_x\times 0$,
where $\R_x\subset \R^2$ is a given axis. 
An admissible tangle is one which has top, bottom, closed or through
strands.  A top strand is a component of the tangle 
that starts at $j^+\in \R_x\times 1$ for some $j$ and
ends at the corresponding $j^-\in \R_x\times 1$, and a bottom strand 
does the same thing at $\R_x\times 0$. A closed strand is a component
$\cong S^1$ in the interior of $\R^2\times [0,1]$. A through strand 
is a pair of components where one component starts at $j^+\in \R_x\times 1$
and ends in $k^{\pm}\in \R_x\times 1$ and the other starts at 
$j^-\in \R_x\times 1$ and ends in $k^{\mp}\in \R_x\times 1$ for some 
$k$ and $j$.

We depict an admissible framed tangle by a generic projection, subject
to the second and third Reidemeister move as well as the usual moves for 
maxima and minima. We will assume the  framing to be in the plane of
projection. $2\pi$-twists in the framing along a strand are depicted
by full or empty blobs as follows:
\picbox{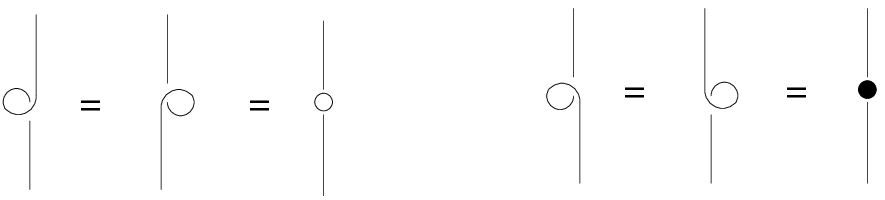}

The admissible tangles are subject to the following relations generalizing 
Kirby's calculus of links \cite{Kir78}:
\ben
\item A Hopf link that is isolated from the rest of the diagram with
one component 0-framed and the other either 1- or 0-framed can be added
or removed from a diagram:
\picbox{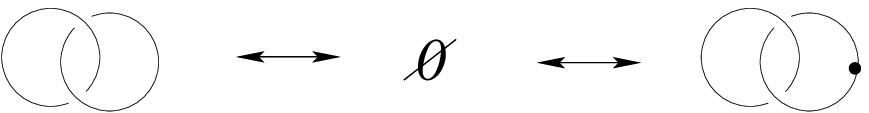}
\item Any strand $R_1$ can be slid over a closed component $R_2$
by a 2-handle slide. This means that $R_1$ is replaced by a connected
sum $R_1\# R_2'$, where $R_2'$ is a push-off of $R_2$ along its framing.  
\picbox{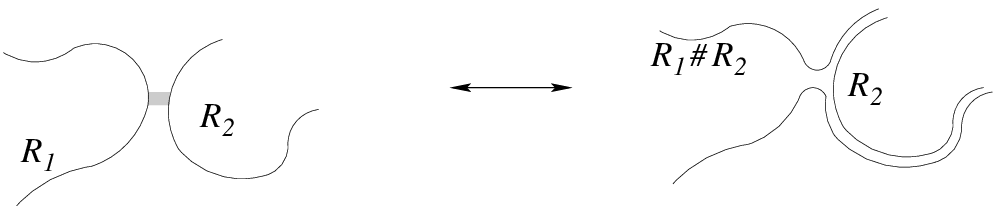}
\item The boundary move, given by introducing two additional components
in a vicinity of points $\{j^+,j^-\}$ at the top line. 
One is an arc connecting $j^+$ to $j^-$, another an annulus going through
that arc and, finally, the outgoing strands are connected through that annulus. 
\picbox{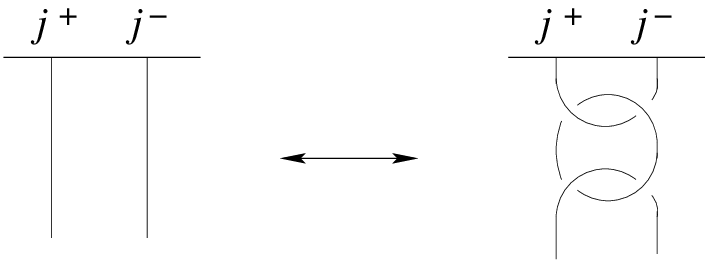}
 \een

Let us record also a few moves that are implied by the above. The first two
can be found in \cite{FR79}. The third follows from the second, the 
2-handle slide  and the boundary move above. 

\ben
\item The Fenn-Rourke Move  in which
a bunch of parallel strands are slid over a 1-framed annulus surrounding them. 
As a  result the group of strands incurs a $2\pi$-twist and a shift in
framing. 
\picbox{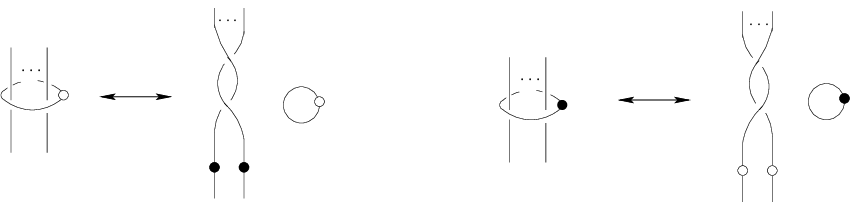}
\item The $\beta$-Move. If 0-framed annulus bounds a disc which intersects
the tangle exactly once with a closed strand then the annulus together
with the closed component can be removed.
\picbox{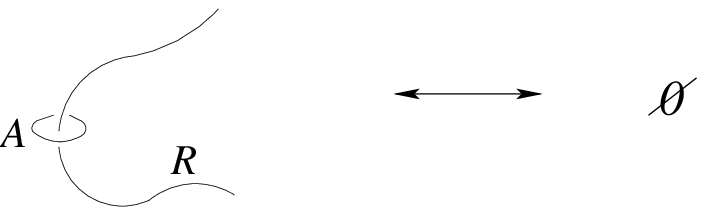}
\item The connecting annulus move. Two different components $R_1$ and $R_2$
of the tangle are linked together by an annulus $A$ as shown. An equivalent 
configuration is the one where $A$ is removed and the two other strands replaced by
$R_1\# R_2$.
\picbox{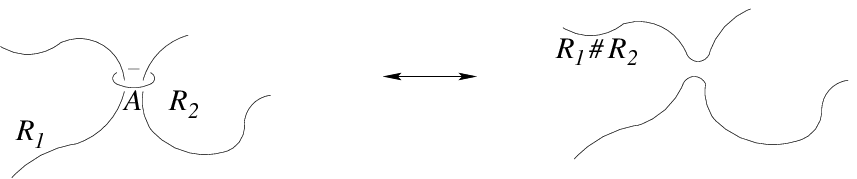}
\een

\paragraph{3.2 Equivalence of braided tensor categories:} 

There is an obvious way in which tangles can be made into a braided 
monoidal category. Two 
tangles $t:k\to l$ and $s:l\to m$ can be composed by stacking 
$t$ on top of $s$ connecting the $2l$ intermediate points with each other. 
It is clear that the composite of admissible tangles is again admissible,
and it is also not hard to prove that this composition factors into
the equivalence classes defined by the moves.  

In addition we have a tensor product $t\otimes u:(k+p)\to(l+q)$
of tangles by putting them side by side into one diagram, that is by
juxtaposing them. Again it follows easily that this operation 
closes in the admissible tangles and factors into the equivalence 
classes. Clearly, this tensor product is strictly associative.

The identity tangle $id:k\to k$ is given by $2k$ parallel vertical strands. 
A braiding $c_{k,l}:(k+l)\to (l+k)$ is given by taking a simple crossing,
as, e.g.,  in \figref{hopfgen}, and replacing one strand by $2l$ parallel
strands and the other by $2k$ parallel strands. All the axioms of a 
braided monoidal category are easily verified. 

\btm[\cite{Ker99}]\label{thm-surg}
There is an isomorphism of braided tensor categories
$$
{\mathfrak Surg}:\;\Tgl\;\stackrel{\cong}{\longrightarrow}\;\Cob\;.
$$
\etm

The assignment of surfaces to numbers is obvious. To produce a cobordism
from a tangle $t:k\to l$ observe first that 
by compactness the tangle must be over some disc
$D^2\subset \R^2$ to which we restrict.  Next we  add
 1-handles to one side of 
$D^2\times[0,1]$ at respective end point pairs and continuing the strand 
through  those handles. Moreover, we bore holes into 
$D^2\times[0,1]\cup\{\mbox{handles}\}\,$ along the strands that start and 
end on the other side.  As a result we obtain a 3-manifold $X$ for which
$\partial X\cong \Sigma_{[k,l]}$. The end points of the tangle have now
disappeared so that we have a link inside of $X$. The desired cobordism
is obtained by performing surgery inside of $X$ along that particular link. 
More details of the construction and the fact that this functor is  well defined 
and an isomorphism are given in  \cite{Ker99}. 

\paragraph{3.3 Generating Tangles:} 

Here we list tangles in \Tgl which are mapped to the generators
in \Cob from Corollary~\ref{cor-gen}. The mapping class group 
generators are as follows:

\picbox{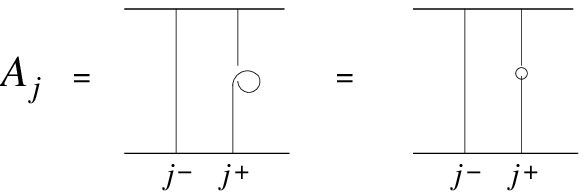}

\picbox{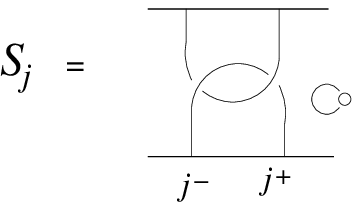}

\picbox{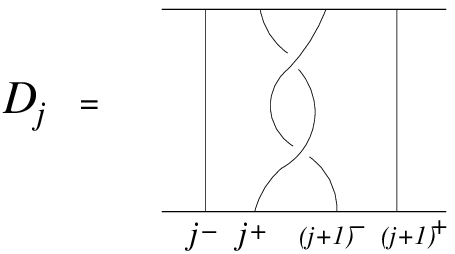}

They are the same as in \cite{MatPol94}. 
The other three generators for handle additions and framing shift are as follows.

\picbox{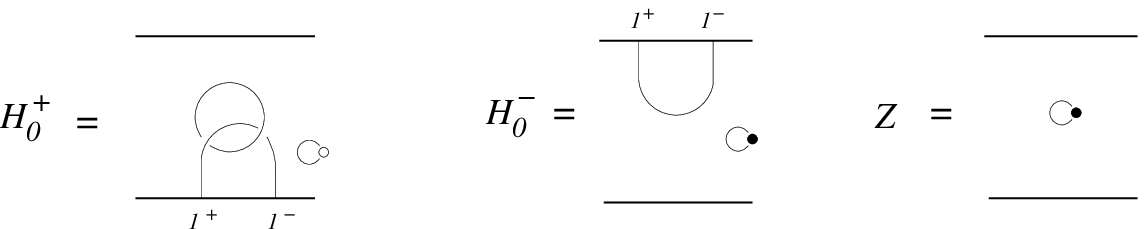}

\head{4.The Category \Alg}\lbl{S4}
The category \Alg is described entirely in algebraic terms.
\begin{definition}
\Alg is the free braided tensor category freely generated by a braided
Hopf algebra object $A$ with two-sided integral and a ribbon element, which 
induces a non-degenerate Hopf paring.
\end{definition}
In particular the set of objects is 
$\{1, A^{\otimes n}\mbox{ with } n\in \N\}$.   
In this section we give a more explicit definition in terms of a
category of planar diagrams. Hence a morphism from $A^{\otimes n}$ to 
$A^{\otimes m}$ in \Alg is represented by a diagram in the strip
$[0,1]\times \R$  with $n$ endpoints at the top end $1\times \R$ and 
$m$  endpoints at the bottom end $0\times \R$. The composite of two
morphisms is given by   stacking the diagrams on top of each other,
and the tensor product of two morphisms by their juxtaposition. 

\paragraph{4.1 The Generators for \Alg :} 
Every morphisms in  \Alg is the product and tensor product of a set of 
generators given by elementary diagrams. They are the  units, $1:1\to A$ 
and $\epsilon:A\to 1$, multiplications,
$\m:A\otimes A\to A$ and $\Delta:A\to A\otimes A$, an invertible  
{\em ribbon element} $v:1\to A$, an invertible antipode,
 $\Gamma^{\pm 1}:A\to A$,  an $S$-invariant 
integral  $\mu:A\to 1$, and an invertible  braid isomorphism 
$c:A\otimes A\to A\otimes A$. The elementary
pictures are the following: 
\picbox{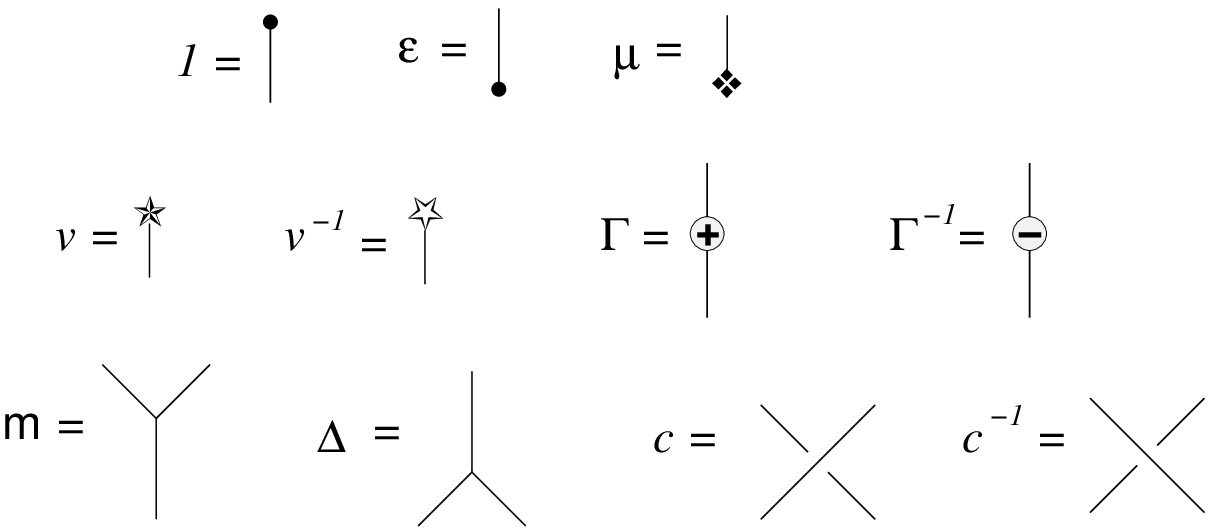}

\paragraph{4.2 The Relations for \Alg :}
The relations for the generating morphisms that define \Alg are mostly
the usual Axioms for braided Hopf algebras. With the conventions as above
we  can  express them as identities between diagrams. 

The first set of identities are those resulting from general isotopies.
This means the Artin braid relations, and the fact that a crossing can
be moved over a fork representing one of the multiplications or over an 
endpoint representing one of elements in $Hom(1,A)$ or $Hom(A,1)$. 

To make $A$ an algebra and coalgebra 
we have to require axioms for associativity and
units, which translate into diagrams as follows. 
\picbox{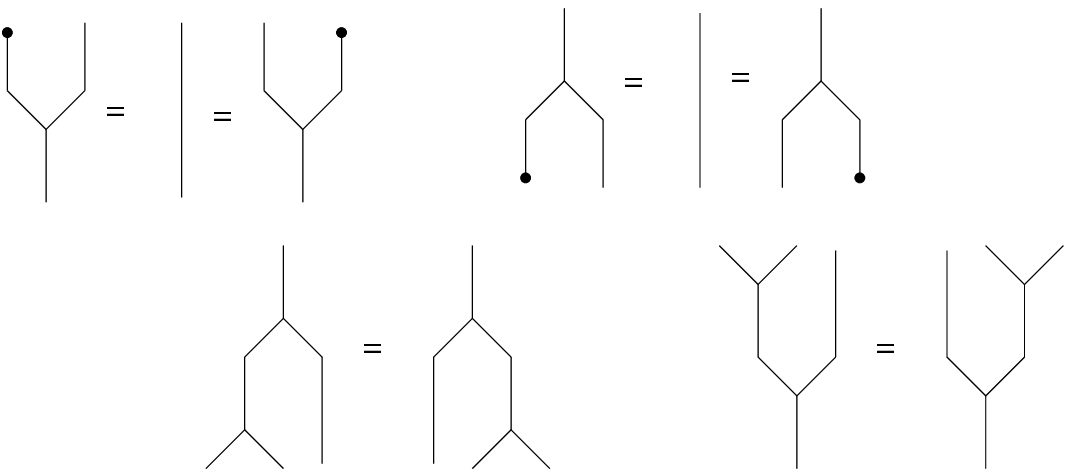}
Next the pictures that make $A$ into a braided bialgebra: 

\picbox{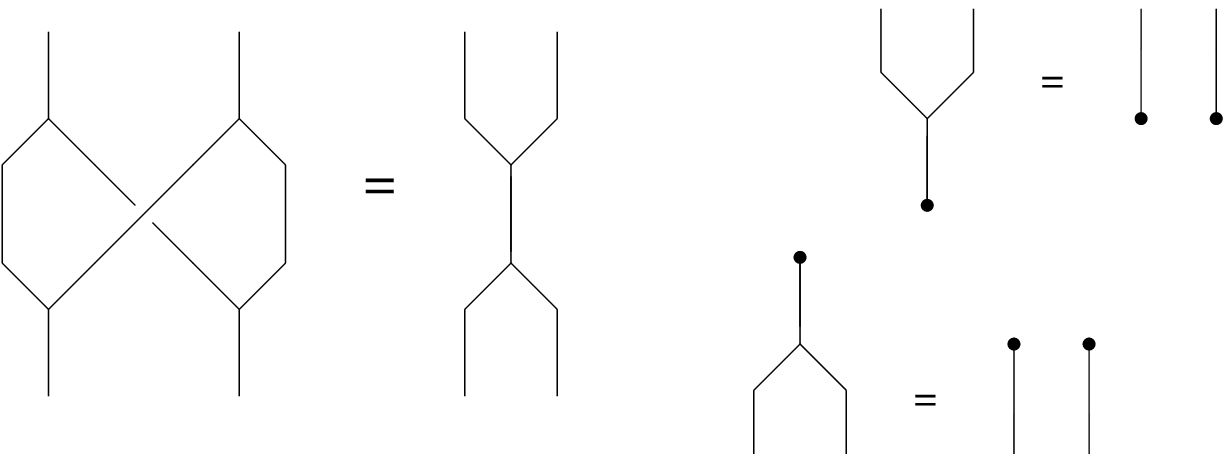}
The axioms for a braided Hopf algebra require the following 
identities for an invertible antipode. 

\picbox{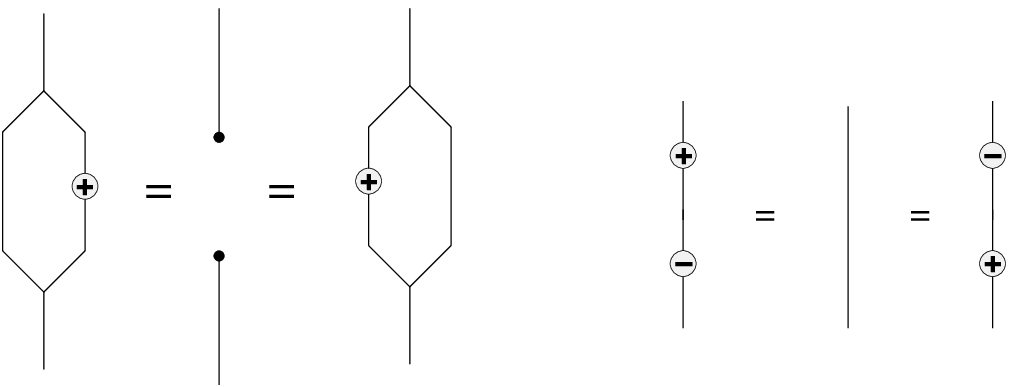}
The defining formula for the right integral and its $S$-invariance  
also have diagrammatic forms. 

\picbox{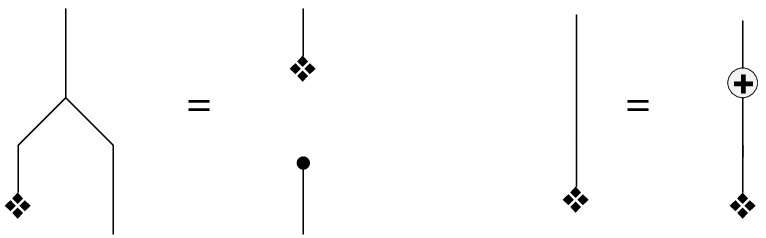}
The ribbon element is firstly required to be central and invertible, 
which is given by the following pictures. 

\picbox{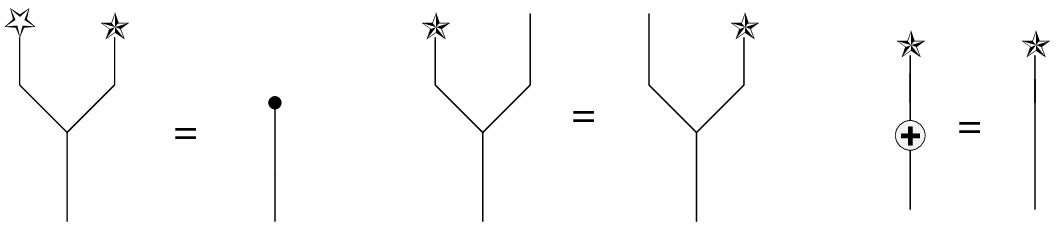}

We denote the operator in the middle picture by $V={\sf m}(v\otimes 1)=
{\sf m}(1\otimes v)$. 
Another property of the ribbon element is that the associated element 
$\omega=V^{-1}\otimes V^{-1}\Delta(v):1\to A\otimes A$ is a Hopf pairing
at least on one side. 

\picbox{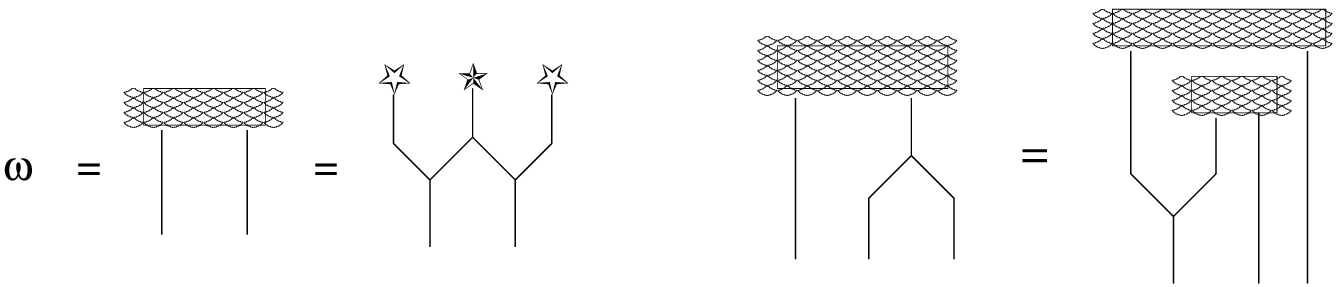}
In this language modularity means that $\omega$ is also non-degenerate.
There are several ways to express non-degeneracy. We will have to require
only a relatively weak version, namely that
$(f\otimes 1) \omega=(g\otimes 1 )\omega$ implies $f=g$. This will imply the
existence of a side-inverse, which is the stronger version. In diagrams
this looks as follows. 

\picbox{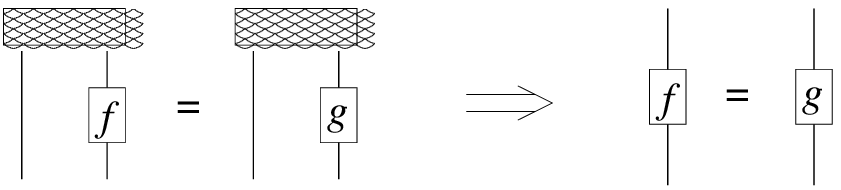}
Finally, we require a number of normalization conditions, which imply
that some morphisms in $Hom(1,1)$ are $1$, meaning their diagrams can
be eliminated. 

\picbox{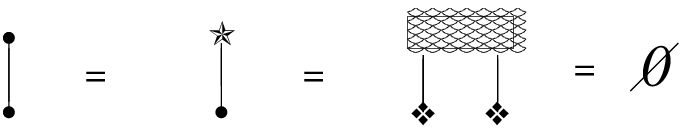}
A  consequence  of this normalization and the previous axioms 
are the following identities.  

\picbox{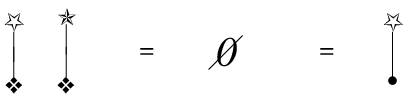}

\paragraph{4.3 Identities for antipode and pairing: } A useful tool
to derive Hopf algebra relation is an algebra structure in
 $Hom(A^{\otimes n},A^{\otimes m})$ given by the convolution product. 
For two morphisms $\alpha,\beta\in Hom(A^{\otimes n},A^{\otimes m})$
we define the product 
$\alpha *\beta=\m^{\otimes m}\circ 
b_m\circ (\alpha\otimes\beta)\circ b_n\circ \Delta^{\otimes n}$,
where $b_n$ are braid morphisms as indicated in the following diagram. 

\picbox{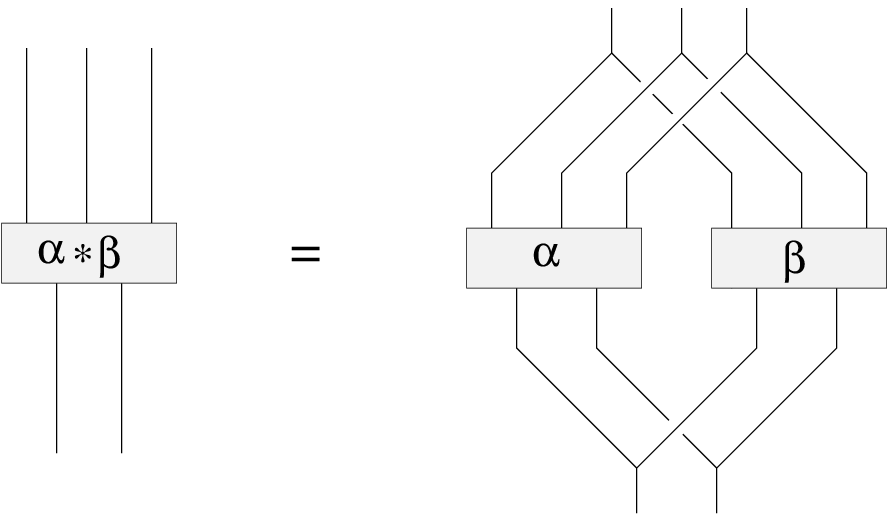}

It is easily checked that the convolution product makes 
$Hom(A^{\otimes n},A^{\otimes m})$ into an associative algebra with
unit $I=1^{\otimes m}\epsilon^{\otimes n}$. Note that the antipode axiom 
in (\ref{fig-hopfanti}) means that $\Gamma$ is the  convolution inverse of the
identity, i.e., $\Gamma*(id)=(id)*\Gamma=I$. Three other natural convolution products
can be found by using inverse braids at the top or bottom half of the diagram.
The first application is a generalization 
Theorem 2.1.4 from \cite{Abe77}:
\begin{lemma}\label{lm-antianti}
\picbox{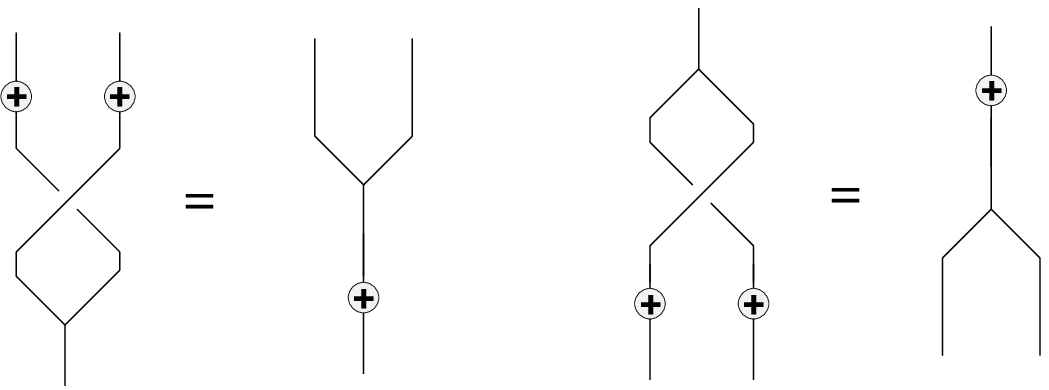}
\end{lemma}
{\em Proof:} As in \cite{Abe77} we consider the operations 
 $L=\m\circ c\circ (\Gamma\otimes \Gamma)$ and $R=\Gamma\circ \m$ on the right and left
hand side of the equation. In the diagram below we compute the
convolution products of these operations in  $Hom(A^{\otimes 2}, A)$ with
$\m$ and find
$L*\m=I=\m*R$. On the left side we use \figref{hopfaxUA} then 
twice \figref{hopfanti} and on the rigth side
\figref{hopfax} and  \figref{hopfanti}.  
\picbox{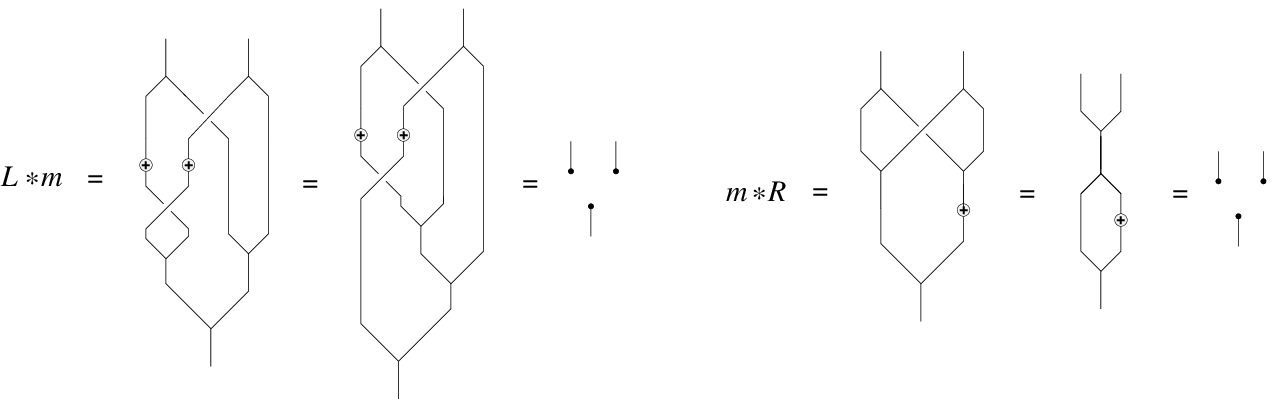} 
This implies $R=L$. The relation for the coproduct follows similarly.

\ep

In addition there are relations for the inverse antipode, in which the 
crossing are exactly opposite. Here the first easy observation, which
follows directly from Lemma~\ref{lm-antianti} and \figref{hopfribb}. 
\blm\label{lm-ribbinv}
The inverse of the ribbon element is central and $\Gamma$-invariant as well.
\elm 
Next let us prove several identities for the form $\omega$.
\blm\label{lm-pairskew} 
The pairing $\omega$ is braided skew and a two-sided Hopf pairing.
In diagrams we have 
\picbox{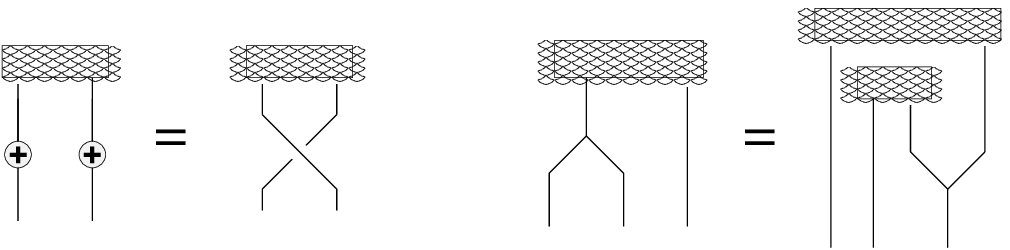}
\elm

{\em Proof:} The skew identity for the antipode is a direct consequence of
Lemma~\ref{lm-antianti}, Lemma~\ref{lm-ribbinv} 
and \figref{hopfribb} applied to \figref{hopfpair}. For the second identity
multiply  $(1\otimes c^{-1})\circ   (c^{-1}\otimes 1)\circ   (1\otimes c^{-1})\circ 
(\Gamma\otimes \Gamma\otimes \Gamma)$ to the two 
sides of the identity in \figref{hopfpair}. 
Applying Lemma~\ref{lm-antianti} and the  skew relation that we just proved
to both sides separately yields the two diagrams above. 

\ep

\blm\label{lm-pairanti} The pairing $\omega^{\dagger}$ defined as $\omega$
but with $v$ and $v^{-1}$ exchanged is given by 
$(\Gamma^{-1}\otimes 1)\circ\omega $.
Moreover, the antipode is self conjugate with respect to $\omega$. This
is summarized in the following picture:
\picbox{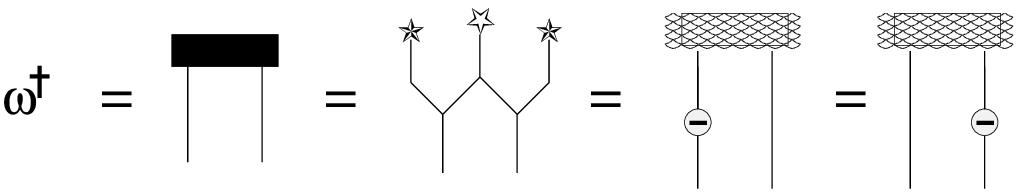}
\elm

{\em Proof:} The first part of the proof is to show that $\omega^{\dagger}$
is a two-sided convolution inverse of $\omega$. That is 
\begin{equation}\label{eq-pairinv}
\omega^{\dagger}*\omega=\omega *\omega^{\dagger}=I.  
\end{equation}
This follows directly from the definition of the convolution product in
$Hom(1,A\otimes A)$, \figref{hopfax}, Lemma~\ref{lm-ribbinv}, 
and \figref{hopfribb}. 

We also find for the convolution product $\omega *((1\otimes \Gamma^{-1})\circ \omega)=I$
by the following diagrammatic calculation. In order we use an isotopy, \figref{hopfpair},
Lemma~\ref{lm-antianti}, \figref{hopfanti}, and \figref{hopfaxUA}. 
\picbox{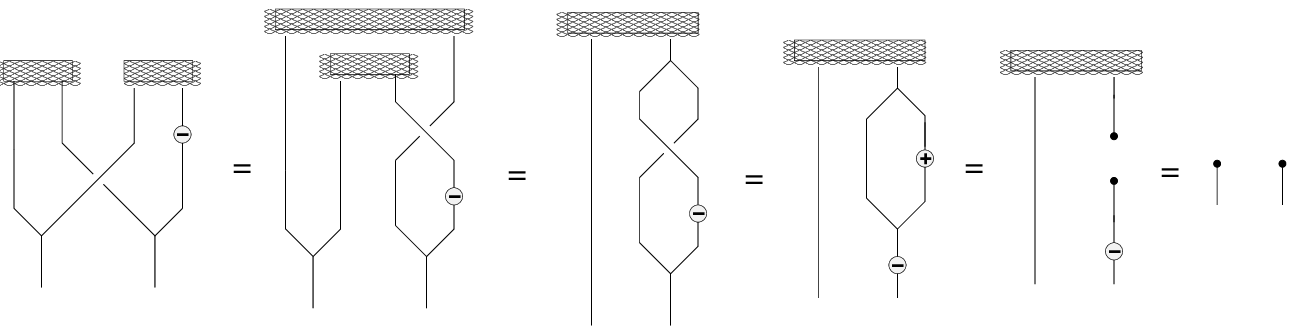}
In the same way we find $((\Gamma^{-1}\otimes 1)\circ \omega)*\omega =I$. Hence 
$\omega^{\dagger}=(1\otimes \Gamma^{-1})\circ \omega=(\Gamma^{-1}\otimes 1)
\circ \omega$.

\ep

\paragraph{4.4. The dual integrals:} We start with a direct consequence 
of Lemma~\ref{lm-antianti}. 
\blm
 $\mu$ is a two-sided integral. That is the reflection of \figref{hopfint}
along a vertical line also holds. 
\elm
 Next we infer the existence of a dual  integral $\lambda:1\to A$. We define it
as $\lambda=(1\otimes \mu)\omega$ and use the following graphical notation.
\picbox{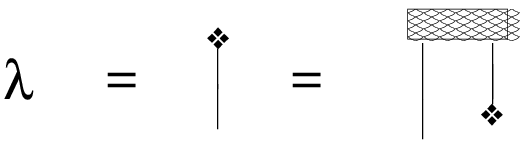}
\blm\label{lm-lambda}  Given that $\omega$ is non-degenerate then 
$\lambda$  is an $\Gamma$-invariant two-sided integral,
and $\lambda=(\mu\otimes 1)\omega$. In pictures 
\picbox{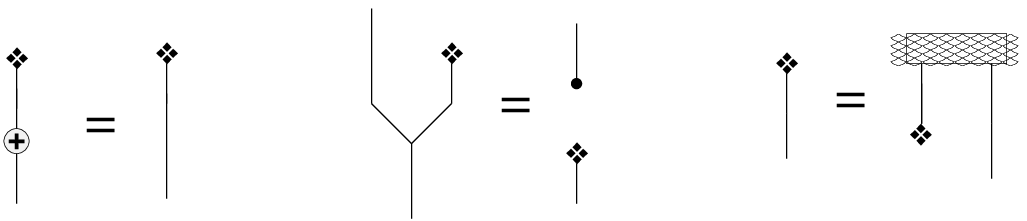}
\elm

{\em Proof:} The fact that $\lambda$ is
 $\Gamma$-invariant follows from self conjugacy of the
antipode, see Lemma~\ref{lm-pairanti}, and 
$\Gamma$-invariance of $\mu$, see \figref{hopfint}. 

The opposite formula $\lambda=(\mu\otimes 1)\omega$ is now a consequence of using the
$S$-invariance of $\mu$ and $\lambda$ at the same time and the skew relation from
Lemma~\ref{lm-pairskew}.
\picbox{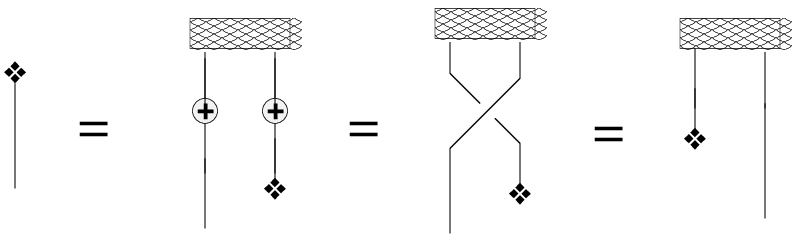}

In the following pictorial calculation we use \figref{hopfpair} and \figref{hopfint}. 
The integral identity for $\lambda$ is then a consequence of the non-degeneracy
condition \figref{hopfnondeg}. 
\picbox{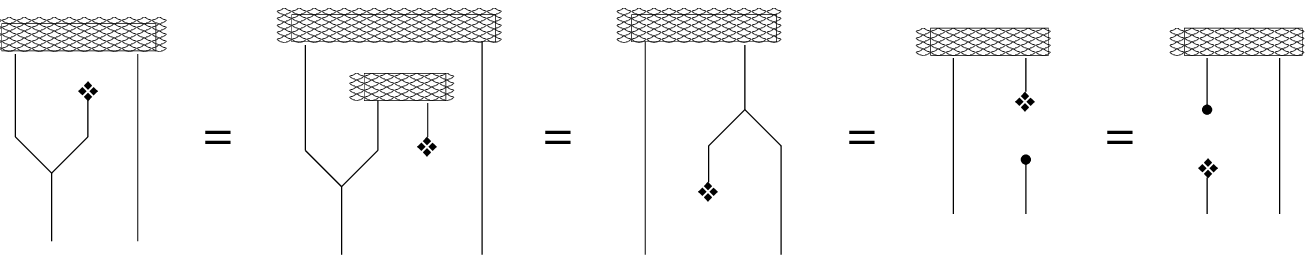}
\ep

Let us next record  an identity that is useful for later calculations. 
It is obtained by considering the identity $(\m\otimes 1)\circ (\Gamma\otimes L)
\circ (\Delta\otimes 1)=
(\m\otimes 1)\circ (\Gamma\otimes R)\circ (\Delta\otimes 1)$,
where $L$ and $R$ are the diagrams on the left and right side of  
the identity on the left side of \figref{hopfax}. 

\picbox{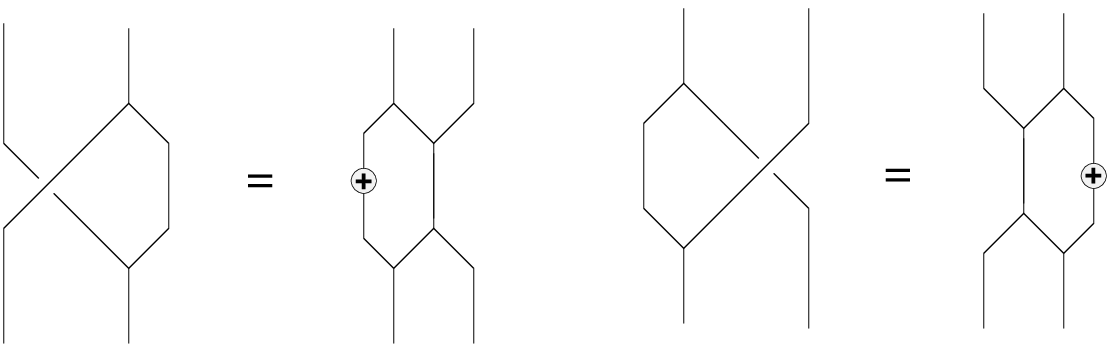}

\paragraph{4.5. Inverse pairings:} A morphisms 
$F\in Hom(1,A^{\otimes 2})$ is the left (right) {\em side-inverse} of some 
 $G\in Hom(A^{\otimes 2},1)$ if $(1\otimes G)\circ (F\otimes 1)=1$
($(G\otimes 1)\circ (1\otimes F)=1$). It follows easily that if a left and 
a right side-inverse exist they must be equal. We determine side-inverses
related to the integrals and Hopf pairings. 

\blm\label{lm-lambdamuinv} Let $\Lambda=\Delta\circ\lambda\in Hom(1,A^{\otimes 2})$ and 
$\Mu=\mu\circ \m \in  Hom(A^{\otimes 2},1)$. 

Then $\Mu(\Gamma\otimes 1)=\Mu(1\otimes \Gamma)$ is the two-sided side-inverse of 
$\Lambda$. Conversely,  $\Mu$  is the two-sided side-inverse of 
$(\Gamma\otimes 1)\Lambda=(1\otimes \Gamma)\Lambda$. In pictures this becomes
\picbox{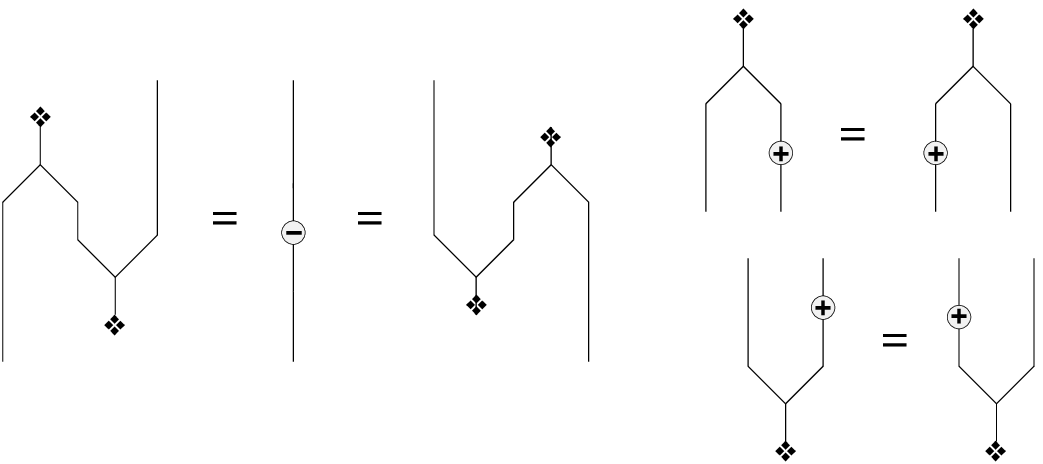}
\elm

{\em Proof:} We multiply to the left identity of \figref{axident} 
$\lambda$ to the upper left end and $\mu$ to the lower right end. 
This is the middle step in the following diagrammatical calculation,
which shows that $(\Gamma\otimes 1)\Lambda$ is a left side-inverse of $\Mu$. 
\picbox{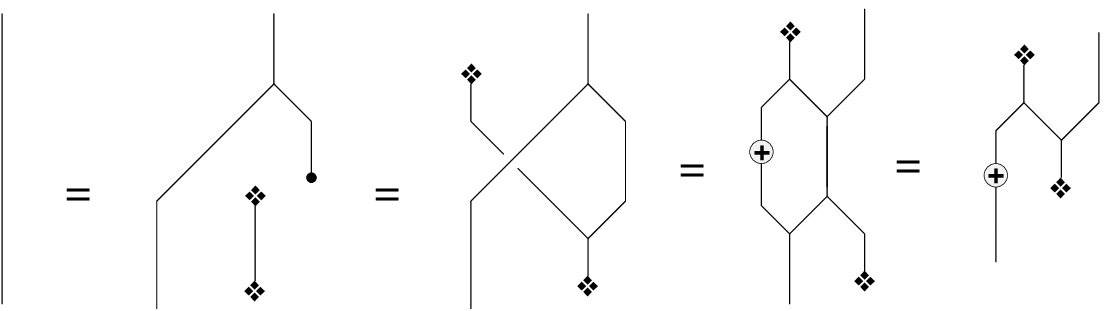}
The other steps are applications of the unit and integral axioms. 
It follows immediately that also $\Lambda$ is 
left side-inverse of $\Mu(1\otimes \Gamma)$. By an analogous calculation
using the second relation in \figref{axident} we find that 
$(1\otimes \Gamma)\Lambda$ is a right side-inverse of $\Mu$, and, further,
that $\Lambda$ is a 
right side-inverse of $\Mu(\Gamma\otimes 1)$. 

\ep

We define a morphism ${\sf S}\in End(A)$ by 
${\sf S}=(1\otimes\Mu)\circ(\omega\otimes 1)$. As a picture we have 
\picbox{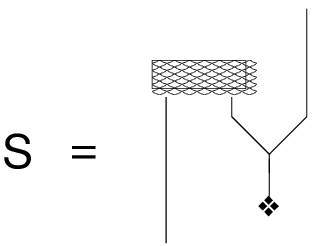} 

\blm\label{lm-SSS} $${\sf S}^2\;=\;\Gamma^{-1}$$\elm

{\em Proof:} This follows from the diagrammatic computation depicted
below. The first picture is the picture of ${\sf S}^2$. The second is 
obtained by applying Lemma~\ref{lm-pairskew}. Next the antipode on the
right is slid through the $\Mu$ and then the $\omega$ pairing using 
Lemma~\ref{lm-pairanti} and Lemma~\ref{lm-lambdamuinv}. An additional
isotopy yields the third picture. The forth follows by an application
of the second part of Lemma~\ref{lm-pairskew}. In the fifth picture we
substitute the definition of $\lambda$ from \figref{deflambda} and use of
Lemma~\ref{lm-antianti}. The last step uses $\Gamma$-invariance of $\lambda$
and Lemma~\ref{lm-lambda} and the pairing identity in
Lemma~\ref{lm-lambdamuinv}.
\picbox{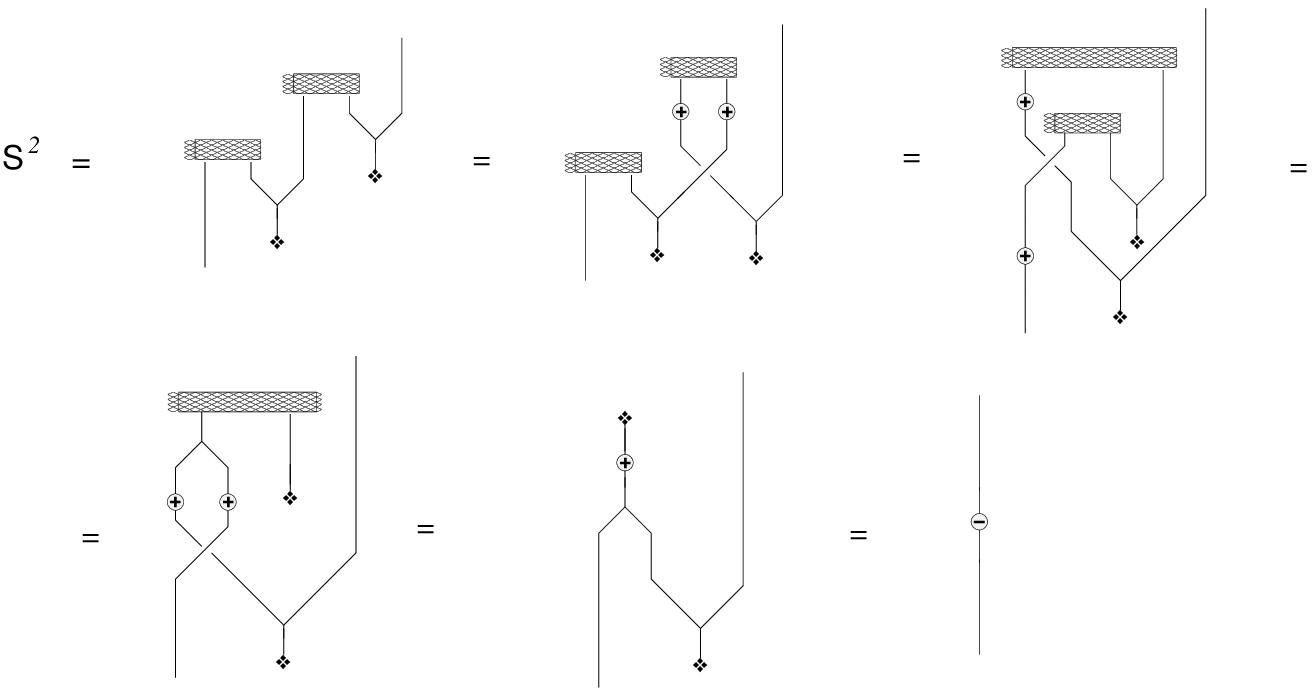}
\ep

Note that Lemma~\ref{lm-SSS} implies the existence of a side-inverse to
$\omega$, which is a slightly stronger property than the non-degeneracy 
from \figref{hopfnondeg}. Conversely, observe that the only place we 
really used the non-degeneracy was in the proof of the fact in Lemma~\ref{lm-lambda}
that $\lambda$ as defined in \figref{lambda} is an integral. If we had required
instead from the start that $\lambda$ is an integral this would have implied
non-degeneracy. Hence the two requirements are interchangeable as axioms!

Let us finally record the following actions of ${\sf S}$ on elements in
$Hom(1,A)$ and $Hom(A,1)$. They are worked out easily using the fact 
that $\omega$ maps units and integrals to each other. The fact that
the inverse ${\sf S}^{-1}$ acts on these elements is a consequence of
their $\Gamma$-invariance and Lemma~\ref{lm-SSS}. 
\begin{equation}\label{eq-Slambda1}
{\sf S}^{\pm 1}\circ \lambda = 1
\end{equation}
\begin{equation}\label{eq-Smueps}
\mu\circ {\sf S}^{\pm 1}  = \epsilon
\end{equation}
\begin{equation}\label{eq-Sribb}
{\sf S}^{\pm 1}\circ v = (\mu\circ v) v^{-1}
\end{equation}

\head{5. A functor from \Alg to \Cob}\lbl{S5}

We construct here the functor ${\mathfrak G}:\Alg\longrightarrow \Cob$.  
This we do by assigning to each generator $g:n\to m$ of $\Alg$  an 
admissible tangle  $T_g:n\to m$ in $\Tgl$. The cobordism is then given by 
${\mathfrak G}(g)={\mathfrak Surg}([T_g])$, where $[T]$ denotes the 
equivalence class  of the tangle $T$ in \Tgl and ${\mathfrak Surg}$ is 
as in Theorem~\ref{thm-surg}.

\paragraph{5.1 Assignment of Generators:}

We begin with the assignments of the product. It is give by the three
component $2\to 1$ tangle written below. 
We can apply the boundary move 
\figref{boundmove} to the left pair
of strands of the left picture 
so that another arc and annulus are introduced. The resulting tangle 
 is symmetric and hence
equivalent to the right picture for ${\sf m}$.  

\picbox{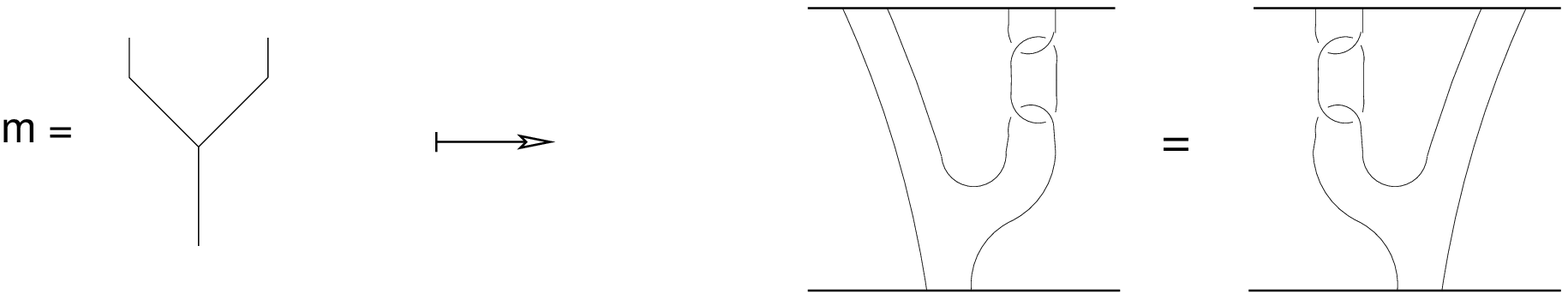}

The coproduct is given by the  next four component $1\to 2$ tangle. 

\picbox{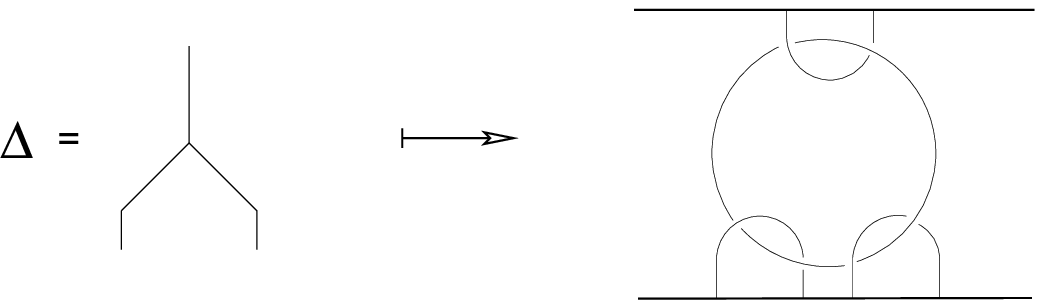}

The antipode $S$ is assigned to the following $1\to 1$ tangle with only
one through pair. Note that the loop can be moved to the other strand by
applying again the boundary move (\figref{boundmove}) on one side, sliding
the loop through and reversing the move \figref{boundmove}. The  
 picture on the right is already the inverse as a braid and so, in 
particular, as a morphism in \Tgl. 

\picbox{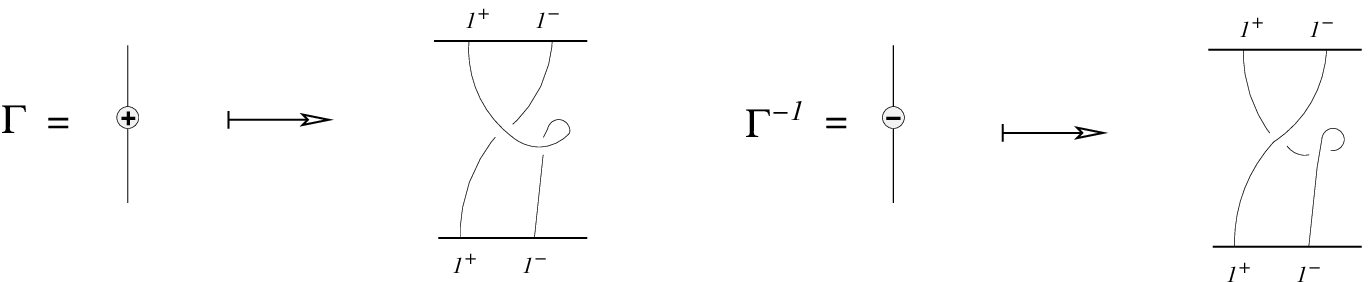}

The generating braid isomorphism $c_{1,1}:1\otimes 1 \to 1\otimes 1$ is 
assigned to the crossing of pairs of strands as follows. Note that,
in our convention, overcrossing strands are mapped to undercrossing ones. 
For the ribbon elements the tangles are given by $0\to 1$ arcs with a
framing loop. 

\picbox{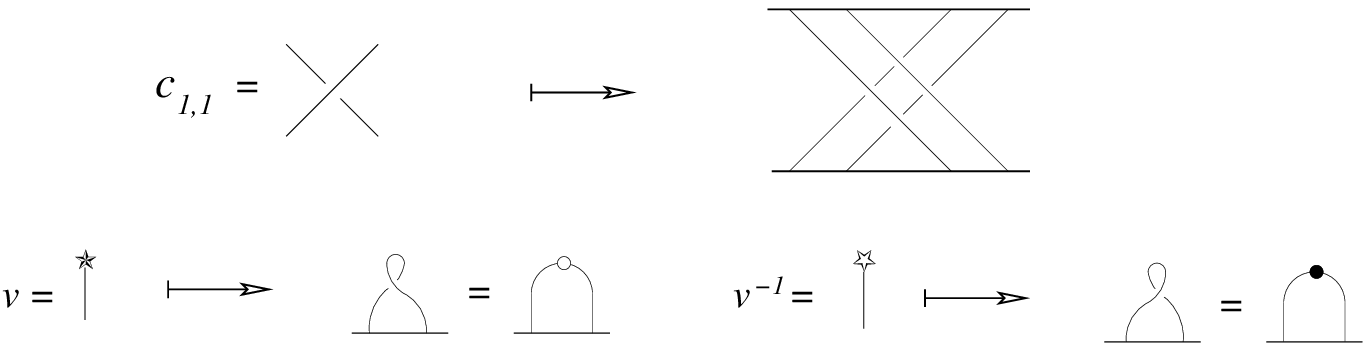}

The unit elements and integrals are also mapped to arcs as in the 
following pictures. We already listed the picture for $\lambda$ here
although it follows from previous assignments. 

\picbox{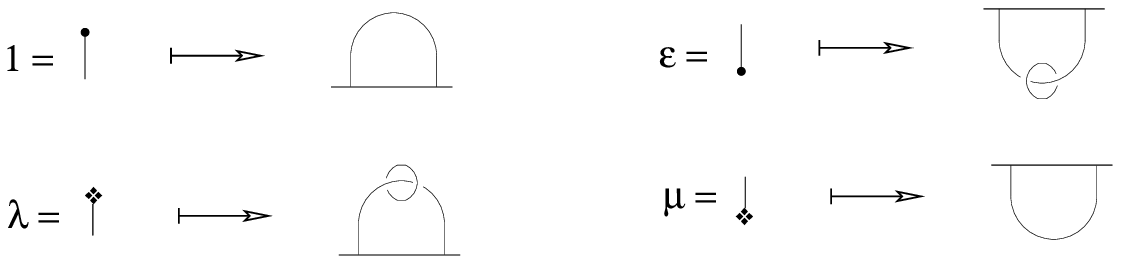}

As cobordisms the units, integrals, and ribbon elements are all homeomorphic
to full tori. 

\blm\label{lm-pairtgl} The pairings are mapped by ${\mathfrak G}$ to the 
cobordism with the following $0\to 2$ tangles. 
\picbox{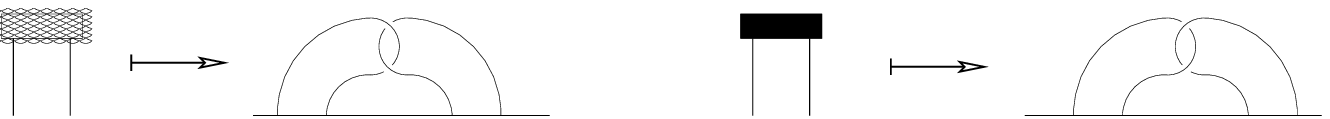}
\elm
{\em Proof:} We start with the tangle associated to 
$V={\sf m}\circ (v\otimes 1):1\to 1$. By composition of the tangles in
\figref{funcprod} and \figref{funcbraid} we obtain the left tangle in
the picture below. The next two pictures follow by a reverse application
of \figref{boundmove} with or without sliding the framing blob through
the annulus. 

\picbox{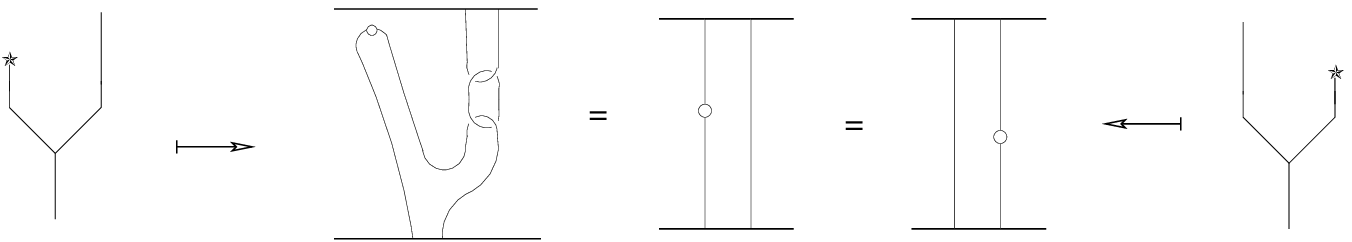}

Next we compute $\Delta\circ v:0\to 2$. The composition of \figref{funccoprod}
and \figref{funcbraid} yields diagram $(a)$. The top annulus is slid off
using  \figref{FRmove} to give $(b)$. Picture $(c)$ is obtained by an
isotopy. We apply \figref{FRmove} again by 
sliding the two arcs at the bottom over the  annulus with the black blob to
give $(d)$. Finally, we obtain $(e)$ by cancelling the isolated annuli as
in \figref{hopflink} and isotopy. 

\picbox{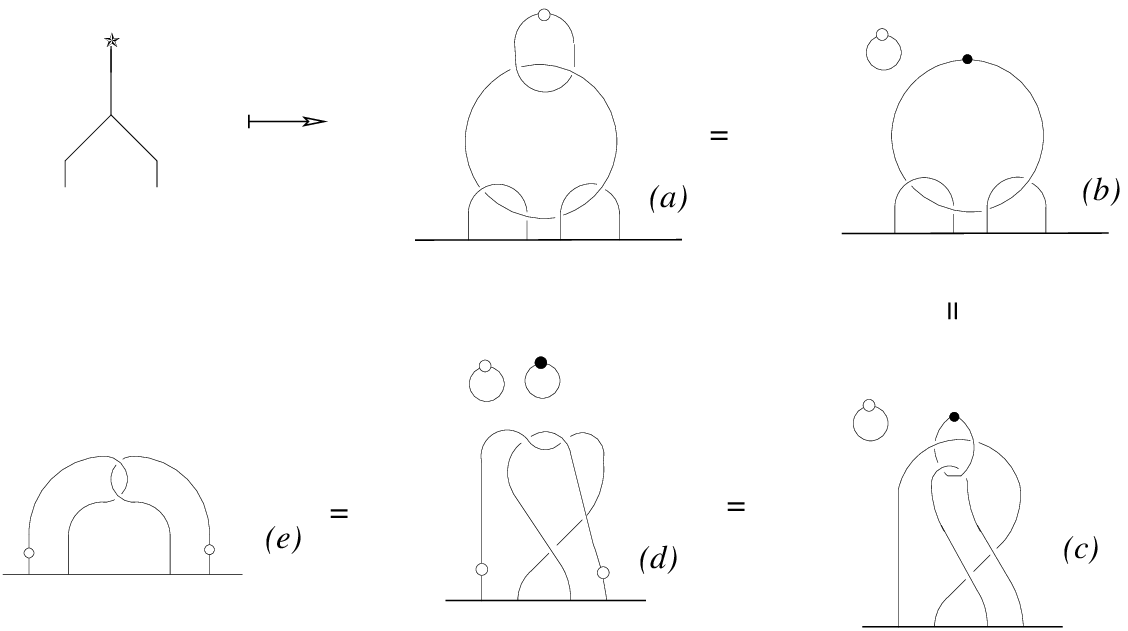}

If we put the expressions from \figref{funcpairPF1} and \figref{funcpairPF2}
together as in the definition \figref{hopfpair} and cancel opposite loops
against each other we find the tangle in \figref{funcpair} for $\omega$.
The proof for $\omega^{\dagger}$ is analogous.

\ep

In the diagrammatic composition below 
we see now immediately that the assignment of $\lambda$ is the one
resulting from Lemma~\ref{lm-pairtgl}. 
\picbox{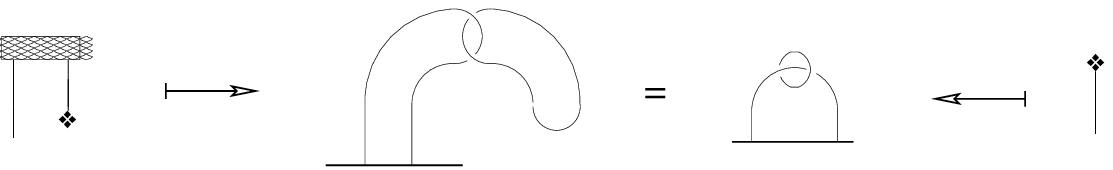}

\paragraph{5.2 Relations:}
  
In order for the assignments of tangles to Hopf algebra generators to
give rise to a functor $\Alg\to \Tgl$ we have to verify that all the
relations in $\Alg$ are also satisfied in $\Tgl$. 

We begin with the unit axioms from \figref{hopfaxUA} for the coproduct. 
The diagram for $(\epsilon\otimes 1)\circ \Delta$ below is obtained by
composition. The pair of annuli $A$ and $R$ is removed with 
the $\beta$-Move
from \figref{betamove}. The resulting three component tangle is equivalent 
to two parallel strands by \figref{boundmove} which represents the identity.
The other three unit axioms follows similarly. 

\picbox{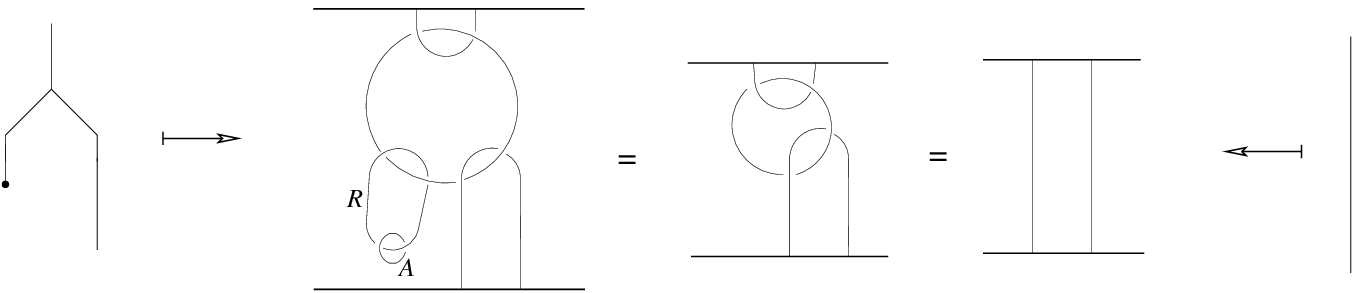}

For the associativity axiom in \figref{hopfaxUA} we pick again only the 
case of the coproduct. The product case is analogous. The pictures for
the two composition $(1\otimes\Delta)\circ\Delta$ and  
$(\Delta\otimes 1)\circ\Delta$ on the left and right. They are equivalent to 
the diagram in the middle by the connecting annulus move from 
\figref{connann}. 

\picbox{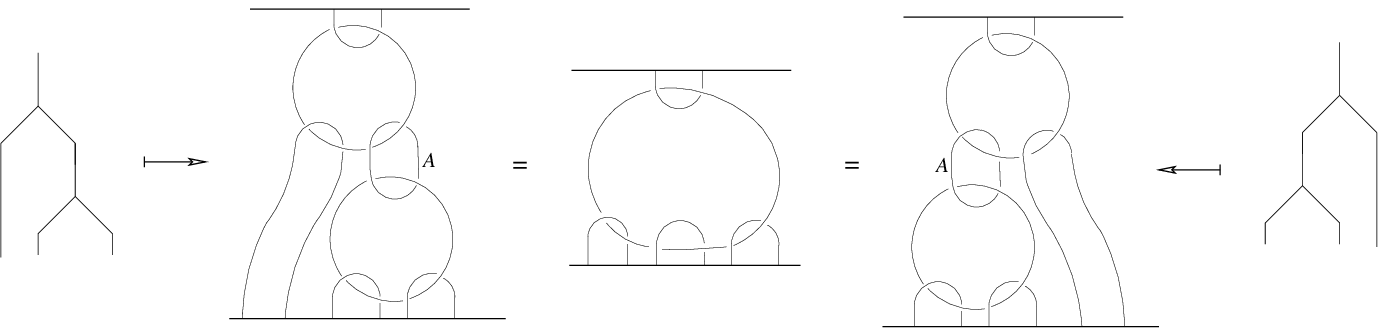}

The bialgebra axiom in \figref{hopfax} is verified next. For 
${\sf m}\otimes{\sf m}\circ (1\otimes c\otimes 1)\circ \Delta\otimes\Delta$
we obtain the tangle in $(a)$ of \figref{Relcompat1} 
by composition. $(b)$ follows by applying
\figref{connann} to each of the annuli in the picture for the ${\sf m}$ part.
An isotopy yields diagram $(c)$ in 
\figref{Relcompat2}. From this we obtain $(d)$ by a 2-handle
slide as in \figref{2slide} of the component labeled $R$ over the component
$S$. The picture in $(e)$ follows again by an isotopy and is precisely
the tangle assigned to $\Delta\circ {\sf m}$.  

\picbox{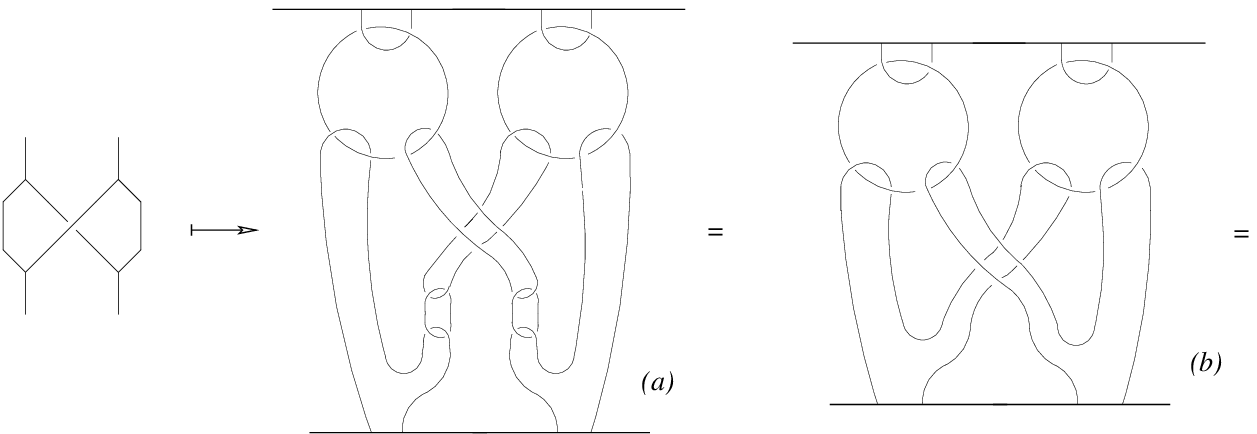}

\picbox{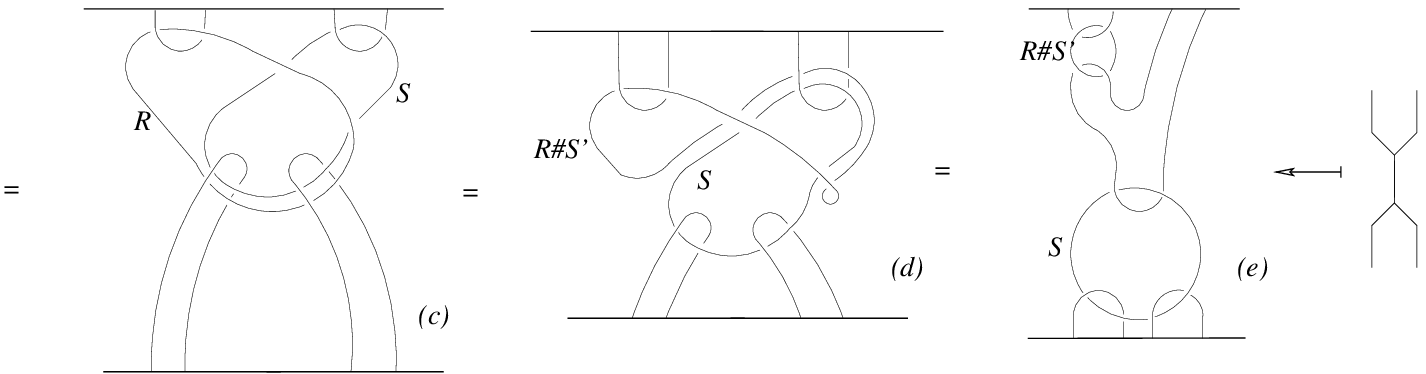}

For the bialgebra axiom with the units in \figref{hopfax} we consider only 
the product case.  The first diagram in \figref{Relbiunit}
is  the composite of the tangles for
$\epsilon\circ {\sf m}$. The second is a result of the boundary move
\figref{boundmove}, and the third follows by applying the $\beta$-Move
from \figref{betamove} to the pair $A$ and $R$. 

\picbox{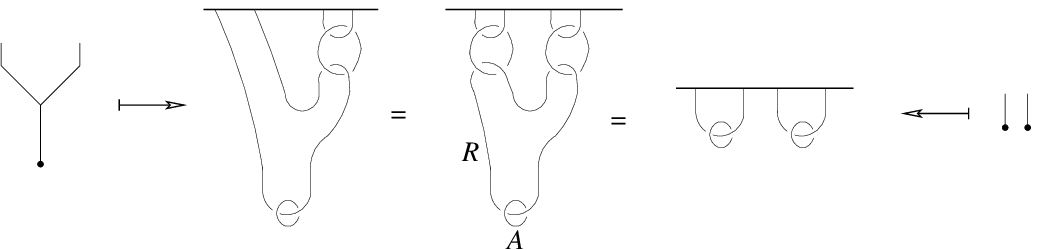}

The antipode axiom is proven in \figref{Relcomp}. To the annulus in the
${\sf m}$-tangle part we apply \figref{boundmove}. The remaining diagrams
follow by isotopies. 

\picbox{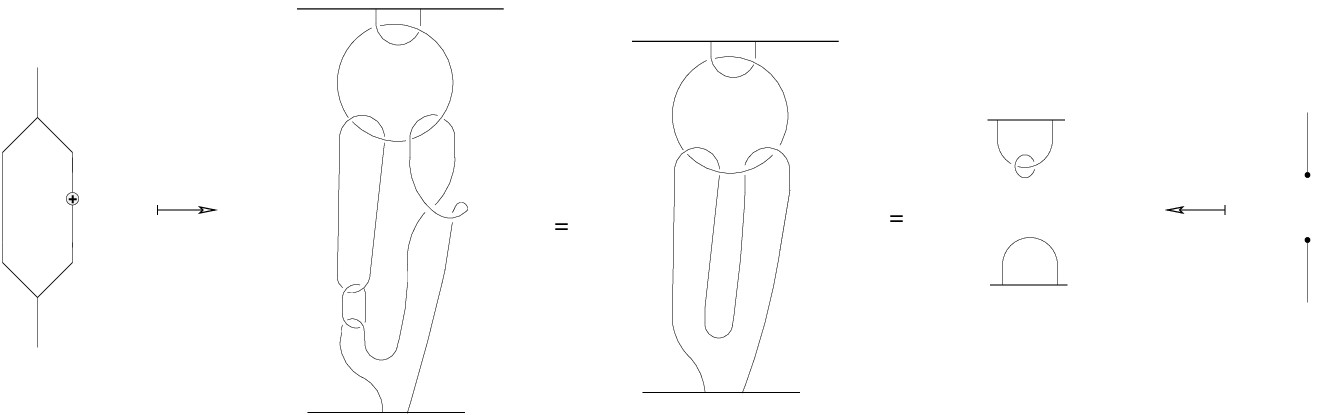}

The axioms for the integrals are verified similarly to those of the units.
One example is given next, which uses again the $\beta$-Move 
\figref{betamove}. 

\picbox{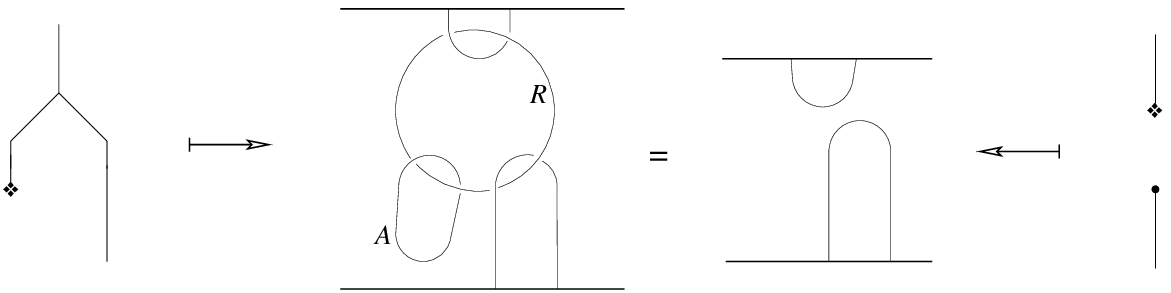}
The pairing axiom from \figref{hopfpair} follows from \figref{Relpair}.
The tangle for the expression in $(1\otimes\Delta)\circ \omega$ is 
depicted in $(a)$. We obtain $(b)$ by application of  \figref{boundmove}
or \figref{connann}. $(c)$ is a result of isotopy and is equivalent to
$(d)$ again by  \figref{connann}. The picture in $(d)$ is also the 
tangle expression for 
$({\sf m}\otimes 1\otimes 1)\circ(1\otimes\omega\otimes 1)\circ \omega$
from \figref{hopfpair}. 

\picbox{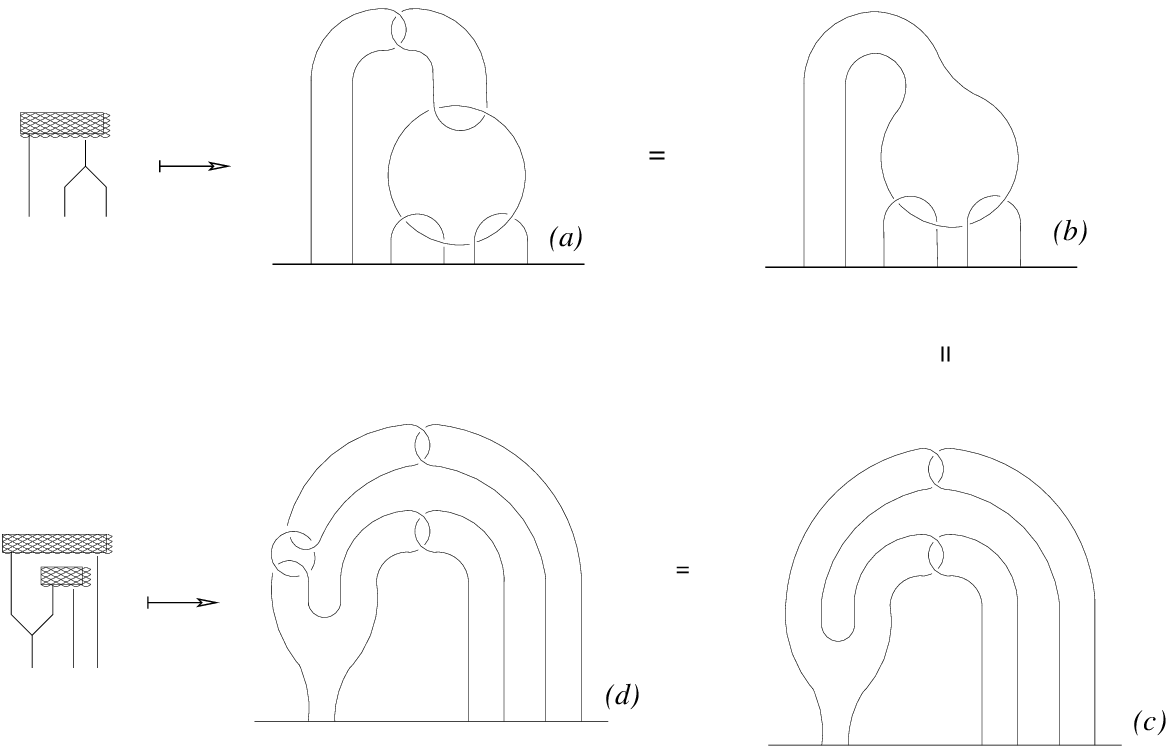}

Non-degeneracy follows easily either by the observation that the tangle
for $\lambda$ has the properties of an integral or by  giving
the side-inverse of tangle for $\omega$ explicitly, namely the reflection
of the tangle for $\omega^{\dagger}$ along the horizontal axis. 

The axioms \figref{hopfribb} for the ribbon element follow already from the
presentations in 
\figref{funcpairPF1} and the first Reidemeister Move for framed tangles,
which allows us to cancel a full with an empty blob on the same strand.

The normalization conditions in \figref{hopfnorm} are a result
of the moves in \figref{hopflink}.  

\begin{cor}\label{cor-algtgl} 
The assignments of tangle classes by representing 
tangles  \figref{funcprod} through 
\figref{funcelem} from  Section~5.1 factor through the relations in \Alg. 

Hence we have a well defined functor\ \  ${\mathfrak X}:\Alg\longrightarrow\Tgl$. 
\end{cor}

Compatibility with the composition and tensor operations are obvious. 
The functor from Theorem~\ref{thm-main} is thus defined as 
\begin{equation}\label{eq-Gfunc}
{\mathfrak G}\;=\;{\mathfrak Surg}\circ{\mathfrak X}\;:\quad\Alg\,
\longrightarrow\,\Cob
\end{equation}

\bigskip

\bigskip

\bigskip

\head{6. From generators  of \Cob to generators of \Alg }\lbl{S6}

Although we cannot construct an inverse functor we will define an 
assignment on the sets of generating morphisms  
\begin{equation}
{\cal W}\;:\;\;Gen[\Cob]\,\longrightarrow\,\Alg
\end{equation}
For the cobordism category we choose them according to 
Corollary~\ref{cor-gen}, that is, 
$Gen[\Cob]=\{H^{\pm}_0, I_{A_1}^{\pm 1}, I_{D_1}^{\pm 1}, I_{S_1}^{\pm 1}, Z^{\pm 1}\}$. 
At least on the set of generators it will be a right inverse for 
${\mathfrak G}$.

\paragraph{6.1 Assignment of Generators and Surjectivity:} 
We give the values of ${\cal W}(g)$ for each $g\in Gen[\Cob]$ and verify 
immediately that for this choice ${\mathfrak G}({\cal W}(g))=g$.

The morphism in \Alg associated to the generator $I_A^{\pm 1}$ is given 
by ${\sf m}\circ (v\otimes 1)$ as depicted in \figref{BF-A}. The tangle 
associated to this is by \figref{funcpairPF1}
 to be exactly the one that represents the
$A_1$-Dehn-twist as in \figref{A-tgl}. 
\picbox{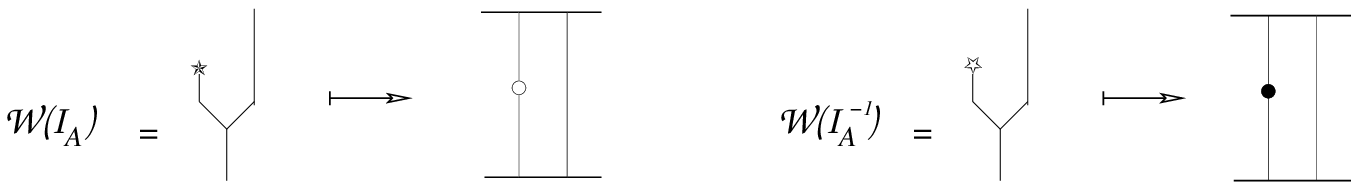}

To the mapping class group generators $I_S$ we associate the morphism
${\sf S}$ from \figref{Sdef} with an additional normalization factor.
The tangle assigned to this composition is depicted in the middle 
of \figref{BF-S}. An application of \figref{boundmove} shows that this
is equivalent to the tangle in \figref{S-tgl}. The morphism in \Alg
that is assigned below to $S^{-1}$ is the inverse of ${\cal W}(I_S)$.
This follows from Lemma~\ref{lm-SSS} and \figref{hopfnormcor}. 

\picbox{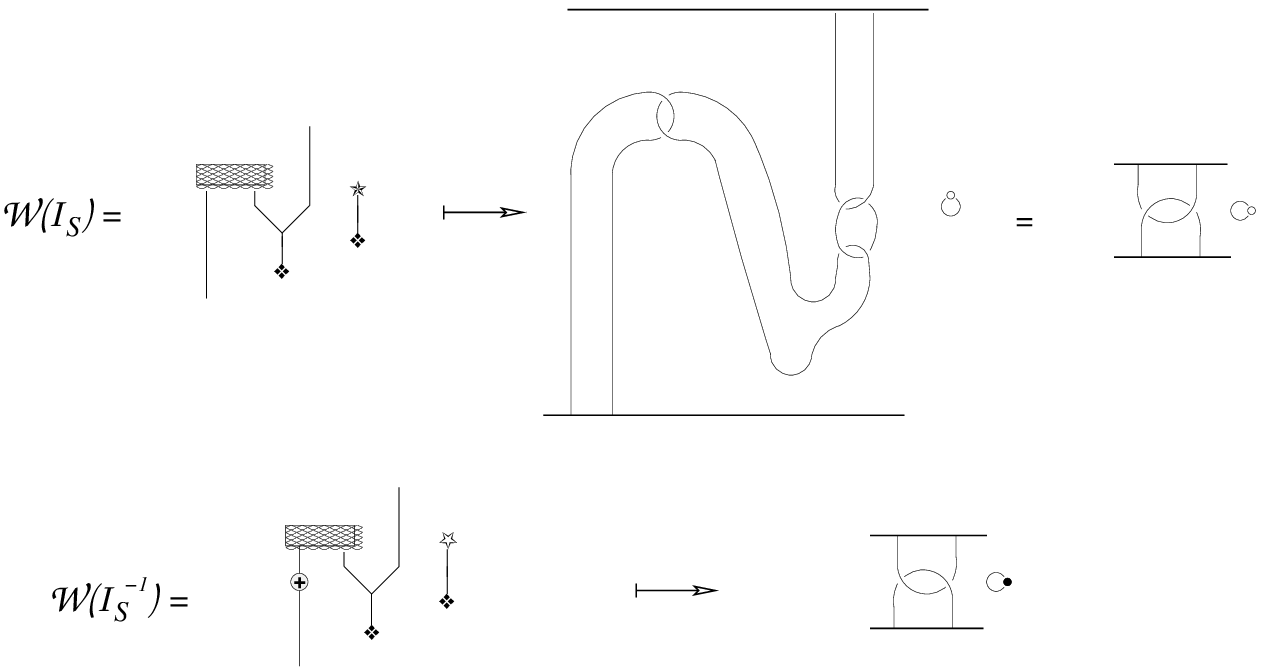}

The last mapping class group generator $I_D$ is mapped to the morphism depicted
in \figref{BF-D}. The next diagram shows the associated tangle for this
morphism. Applying again \figref{connann} or \figref{boundmove} to the
annuli $A$ and $B$ we obtain the tangle \figref{D-tgl}. The morphism we
associate to $I_D^{-1}$ is again the inverse of ${\cal W}(I_D)$ in \Alg.
This is a straightforward exercise using \figref{hopfpair} and then 
\figref{hopfanti}. 

\picbox{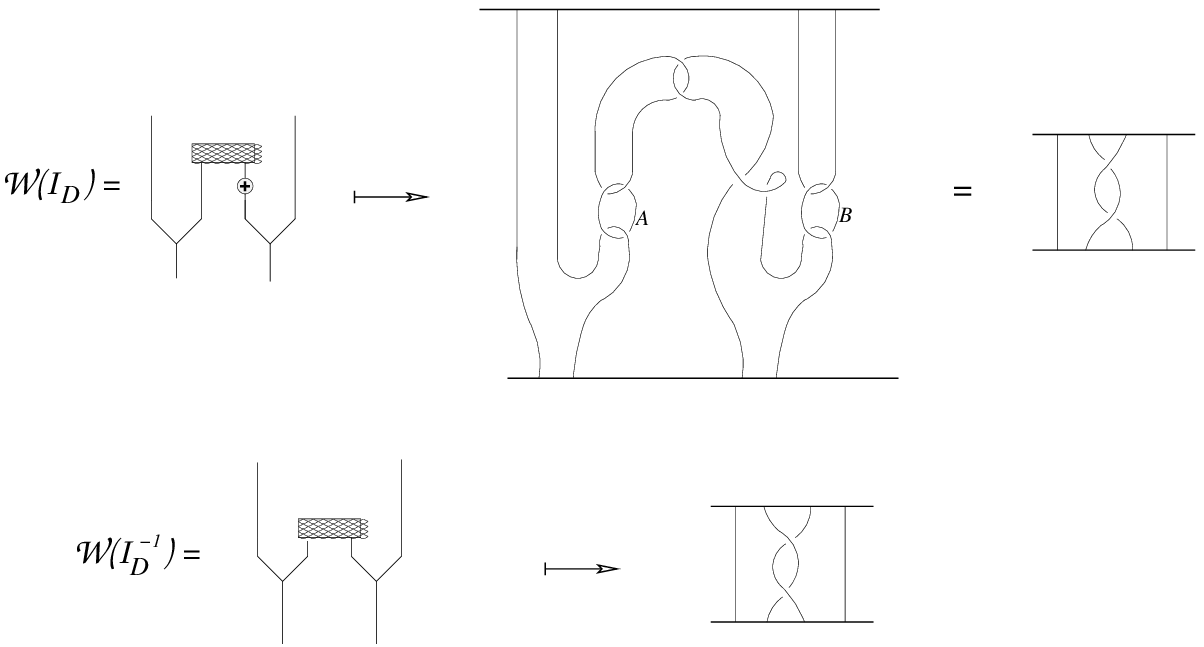}

Finally, we list the assignments for the handle attachments. They are given 
by the integral and normalization pictures. The associated tangles are 
immediately identified with the  pictures in \figref{HG-tgls}. 

\picbox{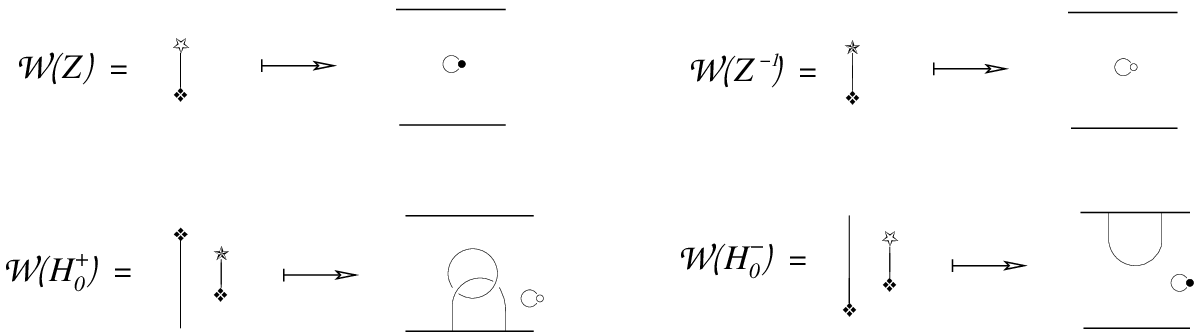}

In summary we found an  assignment ${\cal W}$ vsuch that
\begin{equation}\label{eq-GW=1}
{\mathfrak G}\circ {\cal W}\;\;=\;\;Id\qquad\mbox{on}\quad {Gen[\Cob]}\;.
\end{equation}

Particulary, well definedness of the  functor in (\ref{eq-Gfunc}),
the generators of \Cob in Corollary~\ref{cor-gen}, and the map ${\cal W}$
with the inverse property (\ref{eq-GW=1}) imply now the second part of 
Theorem~\ref{thm-more} as well as Theorem~\ref{thm-main}.

\paragraph{6.2 A Braid Relation:}

A nearby question is whether ${\cal W}$ extends to a functor, and hence
a right inverse to ${\mathfrak G}$, meaning 
${\mathfrak G}\circ {\cal W}$ is the identity on all of \Cob. 
The question is thus, whether all the relations that the generators 
from $Gen[\Cob]$ fulfill in \Cob are also fulfilled by the images in 
\Alg or if we have to introduce additional relations in \Alg. 

Among the set of relations there have to be the relations for the
mapping class group generators. They have been worked explicitly for 
example in \cite{Waj83}. Moreover, we need relations expressing the 
fact that some
mapping class group generators can be extended to the full handles, as
well as relations for Smale-cancellations of handles. Unfortunately, we
do not know of any systematic presentation of \Cob in this way.

Although it might be too optimistic to expect that there are no further
relations in \Alg we demonstrate next that some non-trivial 
relation in $\Gamma_{g,1}$ can indeed be inferred from the relations 
in $\Alg$.  Recall, that the generators 
$A_j$ and $B_j$ of the mapping class group satisfy the braid relation
$A_jB_jA_j=B_jA_jB_j$. Using the definition in (\ref{eq-mcggen}) for $S_j$ 
this relation translates to 
\begin{equation}\label{eq-SArel}
S_j A^{-1}_j S_J\;=\; A_j S_j A_j\;.
\end{equation}

\begin{lemma}\label{lm-ABbrel}
The morphisms ${\cal W}(I_A)$ and ${\cal W}(I_S)$ assigned to the 
generators $A_1$ and $S_1$ fulfill the braid relation (\ref{eq-SArel}).
\end{lemma}

{\em Proof:} The proof is a diagrammatic calculation. The expression 
for the left hand side of (\ref{eq-SArel}) is given in digram $(a)$
of \figref{BraidPF1} below. In $(b)$ we use centrality of the ribbon
elements \figref{hopfribb} and associativity 
to change the order of products. Moreover, we apply  the first
relation in \figref{pairskew}.
The next diagram $(c)$ is the result of an isotopy and the second relation
in \figref{pairskew}. Next, in $(e)$, 
 we use the identity \figref{pairanti} between the regular and 
 the opposite pairing and in $(f)$ the explicit form of $\omega^{\dagger}$
is inserted. In addition we make use of \figref{antianti}. Diagram 
$(g)$ in \figref{BraidPF3} follows then by first cancelling the two
right most ribbon elements, and then applying \figref{hopfint}. We obtain 
an additional factor $\mu\circ v^{-1}$, which we cancel with one of the 
 $\mu\circ v$ factors in the next diagram  $(h)$ using \figref{hopfnormcor}.
Note that it follows from \figref{hopfribb} and \figref{antianti} that
$V$ (multiplication with $v$) commutes with the application of the antipode. 
Hence we can introduce factors $v^{-1}$ and $v$ as indicated in $(h)$. 
The resulting configuration in  the identified in $(i)$ with $\omega$. 
The last diagram $(j)$ follows now by another 
application of \figref{pairanti}. It is readily identified with the 
composite on the right side of (\ref{eq-SArel}).

\picbox{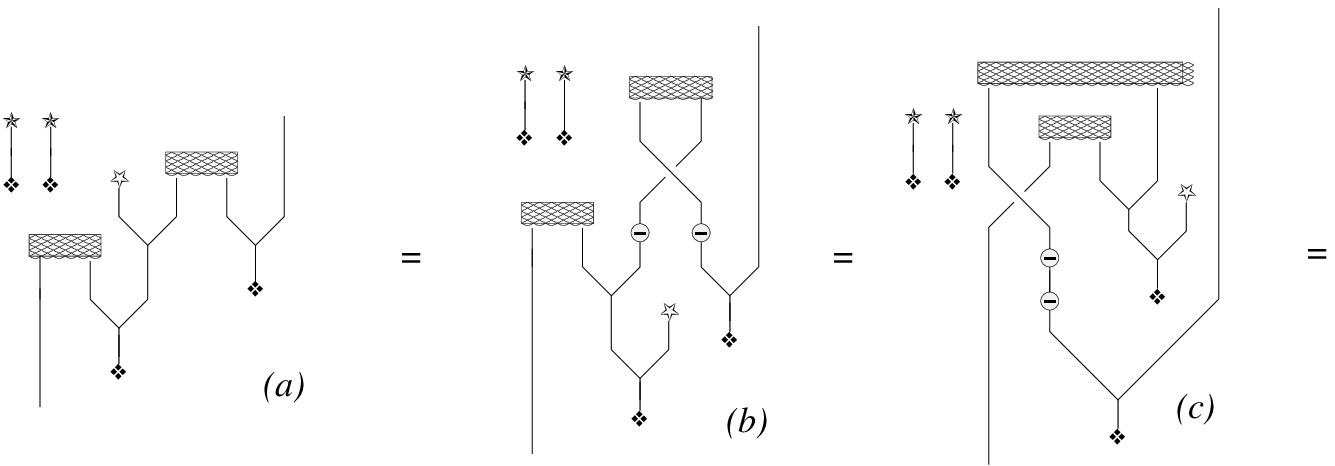}

\picbox{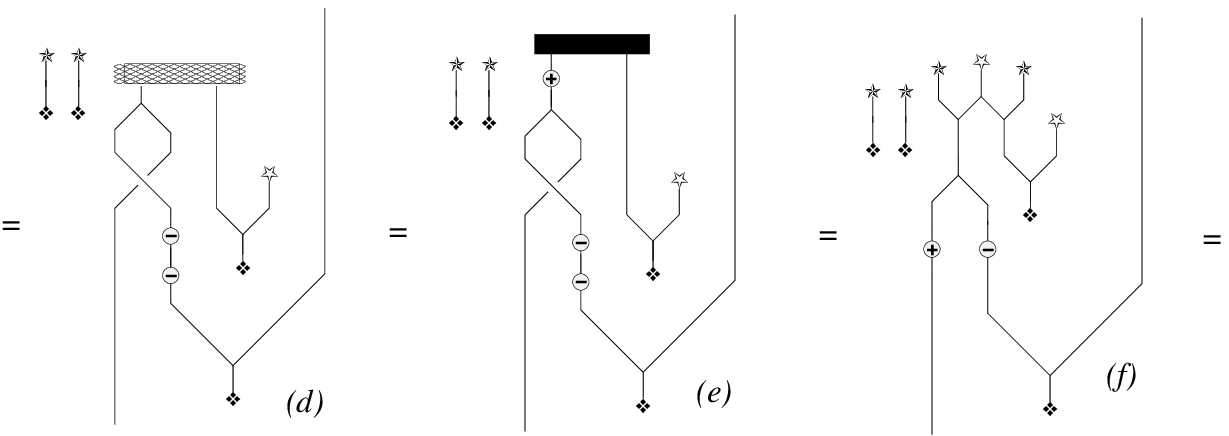}

\picbox{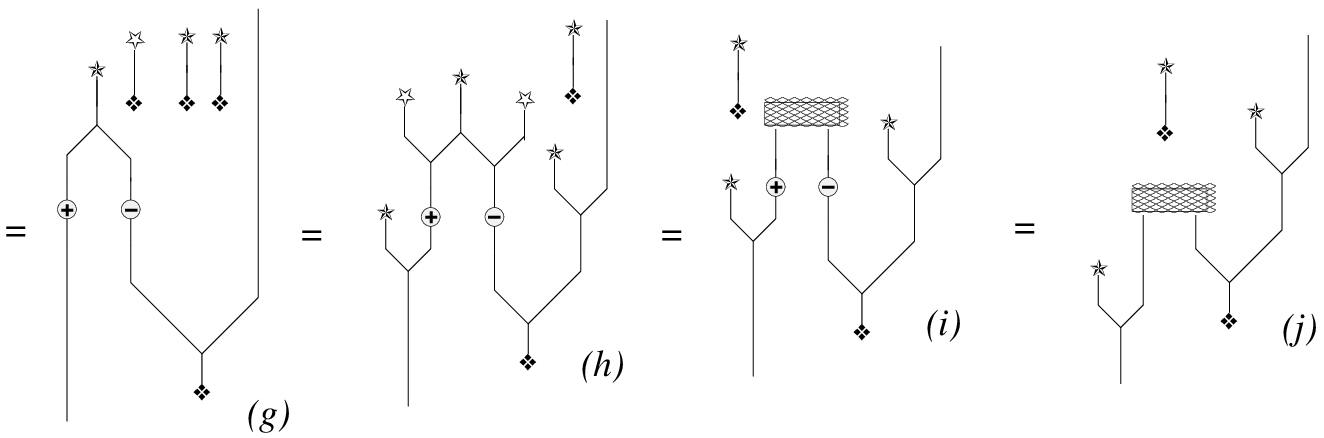}

\ep

The commutation relations for Dehn twist along disjoint are also fulfilled
for obvious reasons. The braid relations between the $B_j$ and $C_j$ are
more difficult to verify and it seems like additional relations have to be
imposed on \Alg. The result is summarized in the first part of Theorem~\ref{thm-more}.

\paragraph{6.3 Heegaard decompositions of generators in \Alg :}

\begin{lemma}\label{lm-multHeeg} The images of the generators under 
${\mathfrak G}$ have the following 
Heegaard decompositions.
$$
{\mathfrak G}(\lambda)=H^+_0\otimes Z
\qquad
{\mathfrak G}(\mu)=H^-_0\otimes Z^{-1}
\qquad
{\mathfrak G}(v)= I_A\circ I_S^{-1} 
$$

$$
{\mathfrak G}(\epsilon)=  H^-_0\circ I_S  
\qquad
{\mathfrak G}(1)=I_S^{-1}\circ H^+_0 
$$

$$
{\mathfrak G}({\sf m})\;=\;(id\otimes H_0^-)\circ I_D 
\circ (id\otimes I_S)  
$$

$$
{\mathfrak G}({\bf \Delta})\;=\;
 (I_S\otimes id)\circ  (I_D^{-1}) \circ (I_S^{-1}\otimes I_S^{-1})\circ 
(id \otimes H_0^+)
$$
\end{lemma}

{\em Proof:} Verification by composition of the associated tangles.

\ep

\begin{propos}
For each of the products of generators of \Cob in 
Lemma~\ref{lm-multHeeg} we have that the 
corresponding product of the images under ${\cal W}$ reproduces the
generators in $\Alg$. 
\end{propos}

{\em Proof:} As a first example we have for the product for 
${\mathfrak G}(\lambda)$ that 
${\cal W}(H^+_0)\circ {\cal W}(Z)=\lambda(\mu\circ v)(\mu\circ v^{-1})$,
the one for ${\mathfrak G}(\mu)$ is similar. For ${\mathfrak G}(1)$
we have ${\cal W}(I_S^{-1})\circ {\cal W}(H^+_0)= 
(\mu\circ v^{-1}){\sf S}^{-1}\circ \lambda(\mu\circ v)=
{\sf S}^{-1}\circ \lambda=1$ by (\ref{eq-Slambda1}). 
The relations for ${\mathfrak G}(1)$ and ${\mathfrak G}(v)$ follow
similarly from (\ref{eq-Slambda1}) and (\ref{eq-Sribb}).

For the multiplication the product $(id\otimes {\cal W}(H_0^-))\circ
{\cal W}(I_D)\circ (id\otimes{\cal W}(I_S) )$ is depicted in 
\figref{REVmult}. The two normalization elements cancel. We also identify 
two pictures for ${\sf S}$ from \figref{Sdef}. Lemma~\ref{lm-SSS} tells
us that we can cancel them together with the antipode. What is left is
${\sf m}$.

\picbox{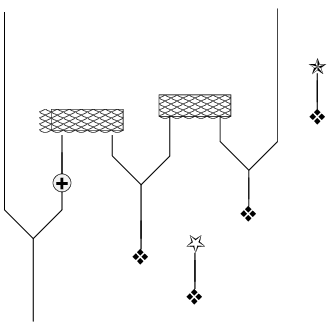}

The situation for the coproduct is more involved. First note that  
${\cal W}(I_S^{-1})\circ{\cal W}( H_0^+)=1$. As a result it suffices 
to show that 
\begin{equation}\label{eq-rrr}
({\cal W}(I_S)^{-1}\otimes id)\Delta = {\cal W}(I_D^{-1})\circ
({\cal W}(I_S)^{-1}\otimes 1)\;. 
\end{equation}
The pictures for the left and right habd side of this equation are
given by the second and third diagram in \figref{REVcoprod1} below. 
Each of them is equivalent to the first and fourth by \figref{pairanti}
and \figref{hopfpair} respectively. 

\picbox{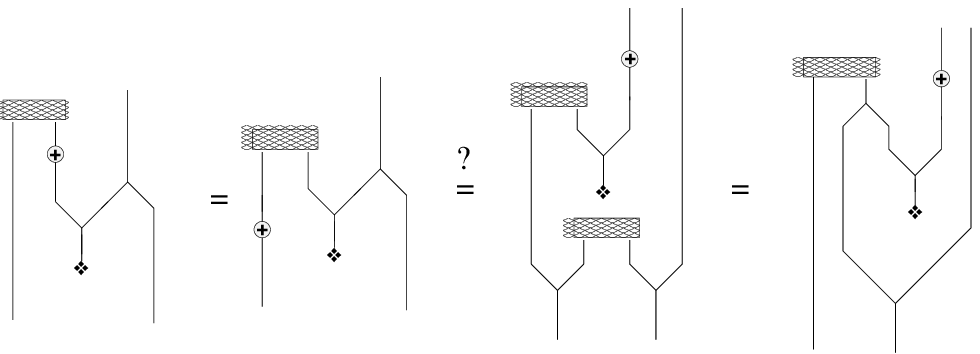}

To check identity of the first and fourth picture in \figref{REVcoprod1}
we can remove the pairing $\omega$ on both sides. Instead let us 
apply the pairing $\Delta\circ\lambda$ to the right most strands on
both sides. We obtain the second and third diagram in  the
next picture \figref{REVcoprod1}. Using Lemma~\ref{lm-lambdamuinv} we
find that the expressions on both sides are equivalent to $\Delta$. 

\picbox{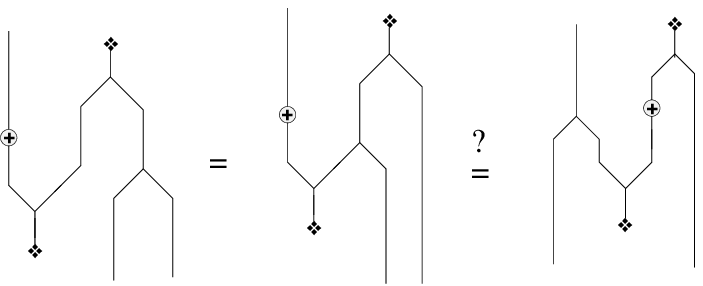}

This proves the equality in (\ref{eq-rrr}).

\ep 

This lemma is useful once we have found a version of \Alg on which
${\cal W}$ extends to a functor
${\mathfrak W}$. We know by (\ref{eq-GW=1}) that
in this case ${\mathfrak G}\circ {\mathfrak W}=id$ but only 
Lemma~\ref{lm-multHeeg} guarantees that 
that ${\mathfrak W}\circ{\mathfrak G} =id$. 

\begin{cor} Let $\overline{\Alg}$ be a quotient of  \Alg, obtained from
the same generators but additional relations. Suppose further
that ${\mathfrak G}$ factors into a functor $\overline{\Alg}\to\Cob$, 
and that ${\cal W}$ extends to a (unique)
functor ${\mathfrak W}$ on $\overline{\Alg}$.
Then ${\mathfrak G}$ and ${\mathfrak W}$ are two-sided inverses and
hence $\overline{\Alg}$ must be isomorphic to $\Cob$.
\end{cor}

This is the same statement as in the last part of Theorem~\ref{thm-more}.

{\sc\small  The Ohio State University, 

Department of Mathematics,

231 West 18th Avenue,

         Columbus, OH 43210, U.S.A. }

 {\em E-mail: }{ \tt kerler@math.ohio-state.edu}
\end{document}